\documentclass[reqno,twoside]{amsart}

\usepackage{amsfonts}
\usepackage{amsmath,amssymb,amsthm}
\usepackage{mathtools}
\usepackage{etoolbox}
\usepackage{color}
\usepackage[english]{babel}
\usepackage{hyperref}
\usepackage{microtype}
\usepackage{fullpage}
\usepackage{floatrow}
\usepackage[utf8]{inputenc}
\usepackage[T1]{fontenc}
\usepackage{lmodern}
\usepackage{textgreek}
\usepackage{comment}

\DisableLigatures{encoding = *, family = * }

\newtheorem{theorem}{Theorem}[section]
\newtheorem{definition}[theorem]{Definition}
\newtheorem{lemma}[theorem]{Lemma}
\newtheorem{proposition}[theorem]{Proposition}
\newtheorem{remark}[theorem]{Remark}

\numberwithin{equation}{section}

\def\RR{{\mathbb{R}}}
\def\NN{{\mathbb{N}}}

\title[]{BOUNDARY OBSERVATION AND CONTROL FOR FRACTIONAL HEAT AND WAVE EQUATIONS}

\author{Umberto Biccari\textsuperscript{\,$\ast$}}  
\address{\textsuperscript{$\ast$}\, Chair of Computational Mathematics, DeustoTech, University of Deusto, Avenida de las Universidades 24, 48007 Bilbao, Basque Country, Spain} 
\email{umberto.biccari@deusto.es}
\thanks{This project has received funding from the European Research Council (ERC) under the European Union’s Horizon 2030 research and innovation programme (grant agreement No. 101096251-CoDeFeL). This material is based upon work supported by the Air Force Office of Scientific Research under award number FA8655-22-1-7012. EZ was partially supported by the Alexander von Humboldt Professorship program; the European Union’s Horizon Europe MSCA project ModConFlex (HORIZON-MSCA-2021-DN-01, project 101073558); the Transregio 154 Project ``Mathematical Modelling, Simulation and Optimization Using the Example of Gas Networks'' of the DFG; and SURE-AI: The Norwegian Centre for Sustainable, Risk-Averse, and Ethical AI, grant 357482, Research Council of Norway. UB and EZ were partially supported by the Grant PID2023-146872OB-I00-DyCMaMod of MICIU (Spain) and by the COST Actions CA24122-Multiscale Stochastics, Patterns, and Analysis of Combinatorial Environments and CA24136-Interactions between Control Theory and Machine Learning. MW is partially supported by the US Army Research Office (ARO) under Award NO: W911NF-20-1-0115} 

\author{Mahamadi Warma\textsuperscript{\,$\dagger$}}  
\address{\textsuperscript{$\dagger$}\,Department of Mathematical Sciences and the Center for Mathematics and Artificial Intelligence, George Mason University.  Fairfax VA 22030 (USA)}
\email{mwarma@gmu.edu}

\author{Enrique Zuazua\textsuperscript{\,$\ddagger$}}
\address{\textsuperscript{$\ddagger$}\, [1] Chair for Dynamics, Control, Machine Learning, and Numerics (Alexander von Humboldt-Professorship), Department of Mathematics, Friedrich-Alexander-Universit\"at Erlangen-N\"urnberg, 91058 Erlangen, Germany.}
\address{[2] Chair of Computational Mathematics, DeustoTech, University of Deusto, Avenida de las Universidades 24, 48007 Bilbao, Basque Country, Spain.} 
\address{[3] Departamento de Matem\'aticas, Universidad Aut\'onoma de Madrid, 28049 Madrid, Spain.}
\email{enrique.zuazua@fau.de}


\begin{document}

\begin{abstract}
We establish boundary null controllability of the heat equation driven by the integral fractional Laplacian $(-\Delta)^s$
on a bounded smooth domain for every $T>0$ and every fractional exponent $s\in(1/2,1)$. The control acts through the singular boundary trace naturally associated with the fractional Dirichlet problem. The main ingredient is a frequency-dependent boundary observability inequality for the associated fractional wave equation, obtained by combining multiplier arguments with the fractional Pohozaev identity. In contrast with the classical wave equation, the observability time deteriorates with the spectral cutoff, reflecting the slow propagation of high-frequency fractional waves. We transfer this estimate to the parabolic problem by transmutation and combine the resulting low-frequency controllability estimates with the high-frequency dissipation of the fractional heat semigroup through a Lebeau-Robbiano iteration. The balance between these two mechanisms yields null controllability precisely in the range $s>1/2$.
\end{abstract}

\maketitle

\section{Introduction}

We prove boundary null controllability for a heat equation driven by the integral fractional Laplacian $(-\Delta)^s$ on a bounded smooth domain. 

A distinctive feature of the problem is that we consider the integral realization of the fractional Laplacian. In contrast with its spectral counterpart, its eigenfunctions are not inherited from those of the classical Dirichlet Laplacian, and standard methods available for local parabolic equations cannot be transferred directly. For this reason, our approach relies on the analysis of the associated fractional wave dynamics, for which we establish a frequency-dependent boundary observation estimate that is subsequently transferred to the parabolic setting through transmutation. 

\medskip
\noindent We consider the fractional heat-like system
\begin{equation}\label{heat_large}
	\begin{cases}
		u_t + (-\Delta)^s u = 0 & \mbox{in }\;\Omega\times (0,T),
		\\
		u=0 & \mbox{in }\; (\mathbb R^N\setminus\overline\Omega)\times (0,T),
		\\
		\rho^{1-s}u=f\chi_{\partial\Omega^+} & \mbox{in }\; \partial\Omega\times (0,T),
		\\
		u(\cdot,0)=u_0 &\mbox{in }\;\Omega,
	\end{cases}
\end{equation}
on a bounded regular domain $\Omega\subset\mathbb R^N$, where $\rho(x)=\mathrm{dist}(x,\partial\Omega)$ denotes the distance to the boundary and, for $s\in(0,1)$ and smooth enough functions $u$,
the operator $(-\Delta)^s$ is formally defined by
\begin{align}\label{fracdef}
	(-\Delta)^s u(x):=C_{N,s}\mbox{P.V.}\int_{\mathbb R^N}
	\frac{u(x)-u(z)}{|x-z|^{N+2s}}\;dz,\;\;x\in\mathbb R^N,
\end{align}
with $C_{N,s}$ an explicit normalization constant.

The function $f$ is localized over a subset $\partial\Omega^+\subset\partial\Omega$, and $ \chi_{\partial\Omega^+}$ stands
for the characteristic function of $\partial\Omega^+$, according to the following partition of $\partial\Omega$
\begin{equation}\label{partition}
	\partial\Omega = \partial\Omega^+\cup\partial\Omega^-, \quad\text{ with }\quad \begin{array}{l}
		\partial\Omega^+ := \{x\in\partial\Omega\,:\,x\cdot\nu>0\},
		\\
		\partial\Omega^- := \{x\in\partial\Omega\,:\,x\cdot\nu\leq 0\},
	\end{array}
\end{equation}
where  $\nu$ denotes the outer normal vector at the boundary $\partial\Omega$.

This boundary action is prescribed through the singular trace naturally associated with the fractional Dirichlet problem, namely through the quantity $\rho^{1-s}u$. This is the framework provided by the existing theory of large solutions for fractional equations (see, e.g., \cite{abatangelo2015large,G-JFA,Ro-Se}), and it gives the natural setting in which a genuine boundary control problem can be formulated for the nonlocal heat equation \eqref{heat_large}.

The main question we want to address is whether the boundary datum $f$ can be chosen so as to drive the corresponding solution of \eqref{heat_large} to rest at time $T$, namely
\begin{align*}
	u(\cdot,T)=0 \qquad \mbox{ in } \Omega.
\end{align*}

When this property holds for every initial datum $u_0$ in a given state space, the system is said to be boundary null controllable in time $T$. In this paper, we prove that this property holds for every $u_0\in L^2(\Omega)$ and every $T>0$ in the range $s\in(1/2,1)$. In fact, our main contribution is the following result.

\begin{theorem}\label{ControlThmHeatIntro}
Let $\partial\Omega^+$ be as in \eqref{partition} and $s\in(1/2,1)$. For any time $T>0$ and any initial datum $u_0\in L^2(\Omega)$, there exists a control function $f\in L^2(\partial\Omega^+\times (0,T))$ such that the unique solution $u$ of \eqref{heat_large} satisfies $u(\cdot,T)=0$ a.e. in $\Omega$. Furthermore, there is a constant $C=C(s,\Omega,T)>0$ such that
\begin{align}\label{control_cost}
	\|f\|_{L^2(\partial\Omega^+\times (0,T))}
	\le C \|u_0\|_{L^2(\Omega)}.
\end{align}
\end{theorem}

Theorem \ref{ControlThmHeatIntro} provides a genuinely multidimensional boundary controllability result for the integral fractional Laplacian. To the best of our knowledge, this is the first result of this type available in the literature for this nonlocal operator.

Establishing Theorem \ref{ControlThmHeatIntro} requires, by duality, a suitable boundary observability inequality for the corresponding adjoint problem. For the integral fractional Laplacian, however, the nonlocal character of the operator prevents the direct application of several standard tools for local parabolic equations, including the Carleman-type approaches underlying classical boundary controllability results \cite{fernandez2000cost,fernandez2000null}. We therefore follow a different route, combining frequency-dependent boundary observability for fractional waves obtained by multiplier techniques \cite{komornik1994exact,lions1988controlabilite,lions1988exact}, transmutation to the parabolic setting \cite{ervedoza1observability,ervedoza2011sharp}, and a Lebeau-Robbiano iteration \cite{Le-Ro}.

This intermediate analysis of the fractional wave equation reveals a further characteristic of the nonlocal dynamics: propagation deteriorates at high frequencies, preventing uniform boundary observability in a fixed time interval. Our approach quantifies this loss and shows how, after transmutation to the parabolic problem, it can be compensated by the high-frequency dissipation of the fractional heat semigroup.

The restriction $s\in (1/2,1)$ reflects this balance. Moreover, this range $s\in(1/2,1)$ in which Theorem \ref{ControlThmHeatIntro} holds is expected to be sharp. Indeed, already for the one-dimensional fractional heat equation \textendash\, which is easier to analyze through spectral methods \textendash\,  the null controllability property fails when $s\in(0,1/2]$ because of the lack of sufficient dissipativity of the system (see \cite{biccari2017controllability,micu2006controllability}).

\subsection{Context}

Partial Differential Equations (PDEs) involving non-local operators, such as the fractional Laplacian, arise naturally in many applications involving anomalous diffusion, stochastic processes, image processing, geophysics, and long-range interactions; see for instance \cite{antil2019sobolev,bakunin2008turbulence,de2002size,dubkov2008levy,gilboa2008nonlocal,gorenflo2007continuous,longhi2015fractional,mandelbrot1968fractional,schneider1990stochastic,vazquez2012nonlinear,weiss2018fractional}.

Among the many analytical questions arising for non-local equations, observation and control have attracted increasing attention in recent years. In fact, despite the extensive literature devoted to the analytical properties of the fractional Laplacian and related non-local models (see, for instance, \cite{biccari2018poisson,BWZ1,biccari2018local,caffarelli2007extension,caffarelli2009regularity,leonori2015basic,Ro-Se,servadei2012mountain,servadei2013variational,War}), control results remain scarce. This is mostly because, in contrast with classical PDEs, many of the standard techniques of control theory relying on finite speed of propagation, localized energy estimates, or Carleman inequalities \cite{coron2007control,komornik1994exact,lions1988controlabilite,tucsnak2009observation}, cannot be applied directly to non-local models. This has motivated the development of new approaches, leading to several controllability results for fractional equations that can be summarized as follows.

\begin{itemize}
	\item[1.] \textbf{One}-\textbf{dimensional problems:} $1$-d problems can be handled by  Fourier series techniques and Ingham type inequalities and the corresponding parabolic versions. In this direction, we refer to \cite{biccari2017controllability,BiWaZu} for the interior null controllability of the $1$-d fractional heat equation, and to \cite{ABPWZ,warma2018null} for the exterior null controllability of the $1$-d fractional heat equation. The case of the $1$-d fractional wave equation with memory terms has been investigated in \cite{BiWa}.
	
	\item[2.] \textbf{Multi}-\textbf{dimensional problems:} the interior null controllability of multi-d fractional Schr\"odinger equations has been treated in \cite{bic} by using the Pohozaev identity for the fractional Laplace operator established in \cite[Proposition 1.6]{Ro-Se}. The exterior controllability of the multi-d strongly damped wave equation has been studied in \cite{Wa-Za}. Apart from that, the approximate controllability for the fractional heat and wave equation in multi-d has been established in \cite{KW-wave,W-wave,WarC}.
	
	\item[3.] \textbf{Spectral fractional Laplacian:} the null controllability of heat equations involving the spectral Dirichlet fractional Laplacian has been thoroughly investigated in \cite{micu2006controllability,miller2006controllability}. The one-dimensional case was established in \cite{micu2006controllability} by means of spectral techniques, while the multidimensional case was obtained in \cite{miller2006controllability} through the Lebeau-Robbiano strategy. In this setting, null controllability holds if and only if $s\in (1/2,1)$, the critical exponent $s=1/2$ marking the loss of controllability. This operator, however, differs substantially from the integral fractional Laplacian considered in the present paper (see, e.g., \cite{servadei2014spectrum} and the references therein), since its eigenfunctions coincide with those of the classical Dirichlet Laplacian and its eigenvalues are simply their fractional powers. Consequently, the existing observability theory for the classical heat equation can be transferred almost directly to the spectral setting, making it considerably more tractable than the integral case considered here.	
\end{itemize}

Nevertheless, additional fundamental questions remain open. In particular, to the best of our knowledge, null controllability of the fractional heat equation in dimension $N \ge 2$ has not previously been established. Theorem \ref{ControlThmHeatIntro} fills this gap. Its proof also yields a frequency-dependent boundary observability result for the associated fractional wave equation, which is of independent interest and reflects the distinctive high-frequency propagation properties of the nonlocal dynamics.

Finally, we mention that an alternative route to the controllability problem considered here could be based on the abstract resolvent approach developed in \cite{miller2006controllability}. There, null controllability of parabolic equations is deduced from suitable resolvent estimates for the generator of the underlying semigroup. Applying this framework to the present setting would require establishing analogous resolvent estimates for the fractional Laplacian together with the boundary observation operator arising from the fractional Pohozaev identity. To the best of our knowledge, however, such estimates are not currently available. 

\subsection{Strategy of the proof}

The proof of Theorem \ref{ControlThmHeatIntro} combines tools from fractional PDEs and control theory and proceeds in four steps.

\medskip 
\noindent\textbf{Step 1: boundary observability of the fractional wave equation.} We first consider the adjoint wave equation
\begin{equation}\label{adjointWaveIntro}
	\begin{cases}
		p_{tt} + (-\Delta)^sp = 0 & \mbox{in }\;\Omega\times (0,T),
		\\
		p=0 & \mbox{in }\; (\mathbb R^N\setminus\Omega)\times (0,T),
		\\
		p(\cdot,0)=p_0,\;\;p_t(\cdot,0)=p_1 & \mbox{in }\;\Omega.
	\end{cases}
\end{equation}

For all $s\in(0,1)$ and $J\in\NN$ fixed, we establish a boundary observability estimate for solutions of \eqref{adjointWaveIntro} involving only the first $J$ eigenfunctions of the Dirichlet fractional Laplacian. Combining the fractional Pohozaev identity of \cite{Ro-Se} with multiplier techniques \cite{komornik1994exact, lions1988controlabilite, lions1988exact}, we show that observability holds for sufficiently large times $T>T_0(J)$, where $T_0(J)=C\lambda_J^{\gamma}$ and the exponent $\gamma=\gamma(s)>0$ depends on the fractional order (see \eqref{gammarange}). In contrast with the classical wave equation, the observability time therefore increases with the spectral cutoff, reflecting the deterioration of propagation at high frequencies.

\medskip
\noindent\textbf{Step 2: boundary observability of the fractional heat equation.} We transfer the frequency-dependent wave observability estimate to the adjoint fractional heat equation
\begin{equation}\label{adjointHeatIntro}
	\begin{cases}
		-v_t + (-\Delta)^sv = 0 & \mbox{in }\;\Omega\times (0,T),
		\\
		v=0 & \mbox{in }\; (\mathbb R^N\setminus\Omega)\times (0,T),
		\\
		v(\cdot,T)=v_T & \mbox{in }\;\Omega.
	\end{cases}
\end{equation}
by adapting the transmutation techniques of \cite{ervedoza2011sharp,ervedoza1observability} based on Kannai's formula (see \cite{kannai1977off}). This yields, for every $T>0$, a boundary observability inequality for solutions involving again only the first $J$ eigenfunctions, with an observability cost depending exponentially on $T_0(J)$.

\medskip 
\noindent\textbf{Step 3: low-frequency boundary controllability of the fractional heat equation.} By duality, the observability estimate obtained in Step 2 yields boundary null controllability of the projection of the fractional heat equation onto the first $J$ eigenfunctions, together with a quantitative estimate of the corresponding control cost. Because the control acts through the singular fractional boundary trace, this step requires a careful formulation of the controlled system in the framework of large solutions, based on transposition techniques.
	
\medskip 
\noindent\textbf{Step 4: controllability of \eqref{heat_large}.} Finally, we combine the quantitative low-frequency controllability estimate of Step 3 with the dissipative decay of the uncontrolled high-frequency components through a Lebeau-Robbiano iteration \cite{Le-Ro}. Unlike the classical construction, which relies on a spectral inequality for packets of eigenfunctions, here we use the frequency-dependent control estimate obtained in the previous steps. The balance between the growth of this control cost and the high-frequency dissipation closes the iteration for $s>1/2$, yielding Theorem \ref{ControlThmHeatIntro}.

\subsection{Structure of the article}
The remainder of the paper is organized as follows. In Section \ref{pohozaev_sec}, we establish a multiplier identity for low-frequency solutions of the fractional adjoint wave equation \eqref{adjointWaveIntro}. This identity is then used in Section \ref{preliminary_sec} to derive the boundary observability of \eqref{adjointWaveIntro}. In Section \ref{sec6}, we employ transmutation arguments to transfer this result to the parabolic setting, obtaining the low-frequency observability inequality  for \eqref{adjointHeatIntro}. Section \ref{sec7} is devoted to proving the boundary null controllability of \eqref{heat_large}. In Section \ref{sec-con-rem}, we present concluding remarks and outline several open problems. Appendix \ref{appendix} collects auxiliary results and properties of the fractional Laplace operator that are used throughout the paper. Finally, Appendix \ref{appendix_proofs} contains the detailed proofs of several technical results.

\section{The Pohozaev identity for the fractional wave equation}\label{pohozaev_sec}

In this section, we consider the following adjoint problem for the fractional wave equation in the absence of control:
\begin{equation}\label{adjointWave}
	\begin{cases}
		p_{tt} + (-\Delta)^sp = 0 & \mbox{in }\;\Omega\times (0,T),
		\\
		p=0 & \mbox{in }\; (\mathbb R^N\setminus\Omega)\times (0,T),
		\\
		p(\cdot,0)=p_0,\;\;p_t(\cdot,0)=p_1 & \mbox{in }\;\Omega.
	\end{cases}
\end{equation}

Our goal is to derive a suitable multiplier identity for this system, that will serve as a key ingredient in establishing the observability properties of the fractional wave equation.

Note that, at this level, the system being time-reversible, taking the initial data at $t = 0$ or $t = T$ is irrelevant, contrarily to the parabolic setting.

The system \eqref{adjointWave}, governed by the non-local fractional operator $(-\Delta)_D^s$ (see \eqref{Ope}), generates a group of isometries. As a consequence, for every pair of initial data $(p_0,p_1)\in H_0^s(\Omega)\times L^2(\Omega)$, it admits a unique weak solution 
\begin{align*}
	p\in C([0,T];H_0^s(\Omega))\cap C^1([0,T];L^2(\Omega))\cap C^2((0,T);H^{-s}(\Omega))
\end{align*}
that preserves the total energy (see, for instance,  \cite[Chapter 10, Theorem 10.14]{brezis2010functional}). Moreover, the function $p$ admits the following spectral representation:
\begin{align*}
	p(x,t) = \sum_{j\in\NN} p_j(t)\phi_j(x)
\end{align*}
where, for each \( j \in \NN \), the coefficient function $p_j$ satisfies the second-order ordinary differential equation
\begin{align*}
	\begin{cases}
		p_j^{\prime\prime}(t)+ \lambda_j p_j(t) = 0,\qquad t\in (0,T),
		\\
		p_j(0) = a_j, \;\;\; p_j'(0) = b_j.
	\end{cases}
\end{align*}

Here, $(\lambda_j,\phi_j)_{j\in\NN}$ denote the eigenvalues and eigenfunctions of the operator $(-\Delta)_D^s$ (see \eqref{e27}), and $(a_j,b_j)_{j\in\NN}$ are the Fourier coefficients of the initial data $(p_0,p_1)$ with respect to the basis $(\phi_j)_{j\in\NN}$, that is,
\begin{align*}
	p_0(x) = \sum_{j\in\NN} a_j\phi_j(x) \quad \mbox{and} \quad p_1(x) = \sum_{j\in\NN} b_j\phi_j(x).
\end{align*}

Here and in the sequel $H_0^s(\Omega)$ denotes the fractional order Sobolev space consisting of all functions in $H^s(\RR^N)$ vanishing in $\RR^N\setminus\Omega$, while $H^{-s}(\Omega)$ is its dual and $\langle\cdot,\cdot\rangle_{-s,s}$ denotes their duality pairing. We will give a more exhaustive description of these spaces in Appendix \ref{appendix} at the end of this paper. 

\smallskip 
\noindent The energy of solutions of \eqref{adjointWave} is conserved along time.  That is, 
\begin{align}\label{e24}
	E_s(t):=\frac 12\int_{\Omega}|p_t|^2\,dx+\frac 12\int_{\Omega}|\nabla^sp|^2\,dx = E_s(0),
\end{align}
where (see \eqref{nablas}) for $p\in H_0^s(\Omega)$ we have set
\begin{align}\label{e25}
	\int_{\Omega}|\nabla^sp|^2\;dx:=\frac{C_{N,s}}{2}\int_{\mathbb R^N}\int_{\mathbb R^N}\frac{|p(x)-p(y)|^2}{|x-y|^{N+2s}}\,dxdy.
\end{align}

Now we establish the  multiplier identity for solutions of \eqref{adjointWave}. To this end, we focus on a special class of solutions $p\in \mathcal H_J$, with 
\begin{align}\label{e32}
	\mathcal H_J = \text{span}\Big\{\phi_1,\ldots,\phi_J\Big\}, \quad J\in\NN.
\end{align}
These are finite energy solutions of \eqref{adjointWave} involving a finite number of Fourier components, of the form 
\begin{align}\label{solFourier}
	p(x,t) = \sum_{j=1}^J p_j(t)\phi_j(x).
\end{align}

The restriction to solutions $p\in \mathcal H_J$ plays a double role throughout the paper. Besides allowing us later to derive frequency-dependent observability estimates, it also guarantees that the multiplier computations leading to the Pohozaev identity are fully justified. Indeed, for every fixed time, $p\in\mathcal H_J$ inherits the regularity properties established in Appendix \ref{appendix} for the individual Fourier components, and all the terms appearing in the identity are well defined. The detailed verification is carried out in Appendix \ref{appendix_proofs}.

The following result holds and will be a key ingredient in the subsequent analysis.

\begin{proposition}\label{lemmaMult}
Let $s\in(0,1)$ and $p\in\mathcal H_J$ be the solution of \eqref{adjointWave} constituted by a finite number of Fourier components, as defined in \eqref{solFourier}. Then,
\begin{align}\label{e210}
	\frac{\Gamma(1+s)^2}{2}\int_0^T \int_{\partial\Omega}(x\cdot\nu)\left|\frac{p}{\rho^s}\right|^2\,d\sigma dt
	= sTE_s(0)+ \int_{\Omega}p_t\left(x\cdot\nabla p + \frac{N-s}{2}p\right)\,dx\,\bigg|_{t=0}^{t=T}, 
\end{align}
where $\Gamma$ is the Euler-Gamma function.
\end{proposition}

\noindent The proof of Proposition \ref{lemmaMult} is presented in Appendix \ref{appendix_proofs}.

\begin{remark}\label{rem22}
{\em 	
The identity \eqref{e210} will be used in Section \ref{preliminary_sec} to obtain a partial observability inequality of \eqref{adjointWave}. To this end, we require sharp estimates for the cross term	
\begin{align}\label{e213}
	\int_{\Omega}p_t\left(x\cdot\nabla p\right)\,dx.
\end{align}
This term poses a significant challenge, as it involves the full gradient of $p$ and therefore cannot be directly controlled by the energy $E_s$ of the solution.
The analysis of this term employs the fact that $p$ only involves a finite number of spectral components and requires us to distinguish various exponent ranges for $s$.	

\begin{itemize} 
	\item \textbf{Range} $s\in(1/2,1)$. Thanks to Proposition \ref{grubb}(c), we have that $\phi_j\in H_0^1(\Omega)$. Thus, since $p$ is constituted by a finite linear combination of eigenfunctions, $p(\cdot, t) \in H_0^1(\Omega)$ as well, for all $t\in[0,T]$, and		\begin{align}\label{216}
		\left|\int_{\Omega}p_t \left(x\cdot\nabla p\right)\,dx\,\right|< C(\Omega)\|p_t\|_{L^2(\Omega)}\|p\|_{H_0^1(\Omega)}.
	\end{align}
	
	\item \textbf{Case} $s=1/2$. From Proposition \ref{grubb}(b) we have that $\phi_j\in H^{1-\varepsilon}(\Omega)$, for all $\varepsilon >0$. This implies that $x\cdot\nabla p\in H^{-\varepsilon}(\Omega)$ for any $\varepsilon >0$ and, for $\varepsilon\in(0,1/2)$, we have that
	\begin{align}\label{217}
		\left|\int_{\Omega} p_t \left(x\cdot\nabla p\right)\,dx\,\right| = \bigg|\langle p_t, x\cdot\nabla p\rangle_{H^{\varepsilon}(\Omega), H^{-\varepsilon}(\Omega)}\bigg| \le C(s,\Omega,\varepsilon)\|p_t\|_{H^{\varepsilon}(\Omega)}\|\nabla p\|_{H^{-\varepsilon}(\Omega)}.
	\end{align}
	Note that $p_t(\cdot,t) \in H^{\varepsilon}(\Omega)$, since $\varepsilon < 1/2$ and $p_t(\cdot,t)$ is also constituted by a finite linear combination of eigenfunctions.
	
	\item \textbf{Range} $s\in[1/4,1/2)$. Proposition \ref{grubb}(a) ensures that $\phi_j \in H^{2s}(\Omega)$. To interpret \eqref{e213} as a duality map, we then need $x\cdot\nabla \phi_j\in H^{-2s}(\Omega)$, i.e. $\phi_j\in H^{1-2s}(\Omega)$, which is satisfied since $2s \ge 1-2s$. Then
	\begin{align}\label{218}
		\left|\int_{\Omega} p_t \left(x\cdot\nabla p\right)\,dx\,\right| = \bigg|\langle p_t, x\cdot\nabla p\rangle_{H^{2s}(\Omega), H^{-2s}(\Omega)}\bigg| \le C(s,\Omega)\|p_t\|_{H^{2s}(\Omega)}\|\nabla p\|_{H^{-2s}(\Omega)}.
	\end{align}	
	
	\item \textbf{Range} $s\in(0,1/4)$. As shown in Step 2 of the proof of Proposition \ref{lemmaMult} (see Appendix \ref{appendix_proofs}), the best estimate that we can provide for \eqref{e213} is 
	\begin{align*}
		\left|\int_\Omega p_t(x\cdot\nabla p)\,dx\,\right|\leq C(\Omega) \int_\Omega |p_t|\,|\nabla p|\,dx \leq C(\Omega)  \|p_t(\cdot,t)\|_{L^\infty(\Omega)}\|\nabla p(\cdot,t)\|_{L^1(\Omega)}.		
	\end{align*}	
\end{itemize} 
	
This estimate will be fully justified employing sharp regularity estimates on the eigenfunctions of the fractional Laplacian and relying on the fact that the solution $p$ under consideration only involves a finite number of modes.

Note that, in each case, the estimate of the cross term \eqref{e213} involves norms of higher regularity than those appearing in the natural energy of the solutions. As a result, the estimates we derive will depend sensitively on the number $J$ of eigenfunctions included in the expansion of $p$. }

\end{remark}

\begin{remark}[\bf Hidden regularity of the fractional normal derivative]
{\em 	
Using the identity \eqref{e210} and the estimates above on the cross term, we can deduce that for every solution $p$ of \eqref{adjointWave} involving a finite number of Fourier components it holds
\begin{align*}
	\int_0^T\int_{\partial\Omega}(x\cdot\nu)\left|\frac{p}{\rho^s}\right|^2\,d\sigma dt < +\infty.
\end{align*}
This represents a \textit{hidden} regularity property for the normal fractional trace $p/\rho^s$, which cannot be obtained through classical trace theorems. It is analogous to the well-known hidden regularity of the normal derivative for solutions to the classical wave equation (see, e.g., \cite{lions1988controlabilite}). However, in the classical case, this hidden regularity ensures that the normal derivative belongs to $L^2$ for all finite energy solutions, whereas in the fractional model considered here, the result only holds for solutions involving a finite number of Fourier modes. This distinction illustrates an important difference with the classical wave equation and points to an issue that warrants further investigation. Related regularity results for solutions of the associated stationary problem can be found in the work of Grubb (see, e.g., \cite{G-JFA} and its references).
 }
\end{remark}

\section{Boundary observability for the fractional wave equation}\label{preliminary_sec}
In this section, we establish frequency-dependent boundary observability results for the fractional wave equation \eqref{adjointWave}. It is important to note that one cannot expect a uniform boundary observability inequality to hold for all finite energy solutions. In the one-dimensional case, this loss of uniform observability can be directly attributed to the failure of the gap condition for the sequence $ (\sqrt{\lambda_j})_{j\in\NN} $, where $\lambda_j$ are the eigenvalues of the operator $(-\Delta)_D^s$; see, for instance, \cite{kwasnicki2012eigenvalues}. While the argument in higher dimensions remains formal, a similar behavior is expected to hold.  Formally, the group velocity of fractional waves is expected to scale like $|\xi|^{s-1}$, which vanishes as $|\xi| \to +\infty$. This decay in propagation speed at high frequencies presents a main obstruction to the validity of the observability inequality for general solutions.

This is closely related to the previously mentioned observation that the estimates for the cross term in \eqref{e213} involve norms of higher regularity than those associated with the system's natural energy.

For this reason, we focus on the specific class of low-frequency solutions $p\in \mathcal H_J$, for which boundary observability results can be established, but with a dependence on the frequency number $J$. In more detail, we have the following result.

\begin{proposition}\label{prop31}
Fix $s\in(0,1)$ and $J\in\NN$. If $s\in(1/2,1)$, fix in addition $\varepsilon\in(0,s-1/2)$. Set
\begin{equation}\label{e39}
	T_0(J):=C\lambda_J^{\gamma},
\end{equation}
where
\begin{equation*}
	C=
	\begin{cases}
		C(N,s,\Omega), & s\in(0,1/2],\\
		C(N,s,\Omega,\varepsilon), & s\in(1/2,1),
	\end{cases}
\end{equation*}
and
\begin{align}\label{gammarange}
	\gamma=
	\begin{cases}
		\displaystyle \frac{1-s}{1-2\varepsilon}, & \displaystyle\text{if } s\in\left(\frac12,1\right),\\[10pt]
		1, & \displaystyle\text{if } s=\frac12,\\[10pt]
		\displaystyle 1+\frac{\alpha(s)}{2}, & \displaystyle\text{if } s\in\left[\frac14,\frac12\right),\\[10pt]
		\displaystyle \frac Ns+1, & \displaystyle\text{if } s\in\left(0,\frac14\right).
	\end{cases}
\end{align}
with $\alpha(s)$ given by \eqref{eq:alpha}. Let $\partial\Omega^+$ be as in \eqref{partition}. Then, for all $T>T_0(J)$, the boundary observability inequality
\begin{align}\label{e37}
	E_s(0)\le \frac{\Gamma(1+s)^2}{2s(T-T_0(J))}
	\int_0^T\int_{\partial\Omega^+}(x\cdot\nu)\left|\frac{p}{\rho^s}\right|^2\,d\sigma dt
\end{align}
holds for every solution $p\in\mathcal H_J$ of \eqref{adjointWave} defined in \eqref{solFourier}.
\end{proposition}

\noindent The proof of Proposition \ref{prop31} is presented in Appendix \ref{appendix_proofs}.

\begin{remark}\label{obsRem}
{\em
Let us stress that the minimal observability time $T_0(J)$ given in \eqref{e39} depends on the frequency number $J$. 
	
When $s = 1$, i.e., in the case of the classical wave equation, the observability time $T_0$ is independent of the frequency, which aligns with the uniform finite speed of propagation characteristic of the classical wave equation and its well-established observability properties (see \cite{lions1988exact}). In contrast, when $s \in (0,1)$, the observability time $T_0(J)$ diverges as $J \to +\infty$. This behavior is consistent with the earlier observation on the lack of uniform velocity of propagation when $s \in (0,1)$; see, for example, \cite{kwasnicki2012eigenvalues}.}
\end{remark}

\section{Transmutation: observability of the fractional heat equation} \label{sec6}

In \cite{ervedoza1observability,ervedoza2011sharp}, the authors introduced a variant of  the so-called \textit{transmutation techniques}, which allow for the transfer of observability results from the hyperbolic to the parabolic setting. In this section, we demonstrate how these techniques can be applied to the boundary observability inequality \eqref{e37} for the fractional wave equation, leading to an analogous observability result for the fractional heat equation. To this end, we will make use of the following result taken from \cite[Proposition 2.4]{ervedoza1observability} (see also \cite{ervedoza2011sharp}).

\begin{lemma}\label{lemmaKernel}
Let $T>0$ and $L>0$ be given real numbers. Then, for every $\beta>2L^2$, there exists a function $k=k(\zeta,t)$ satisfying 
\begin{align}\label{eq:k}
	\begin{cases}
		k_t + k_{\zeta\zeta} = 0 &\mbox{in }\;(-L,L)\times (0,T),
		\\
		k(\zeta,0) = k(\zeta,T) = 0 &\mbox{in }\;(-L,L),
		\\
		k(0,t) = 0 &\mbox{in }\; (0,T),
	\end{cases}
\end{align}
and 
\begin{align}\label{k_id}
	\displaystyle k_\zeta (0,t) = \exp\left[-\beta\left(\frac 1t + \frac{1}{T-t}\right)\right], \quad t\in(0,T).
\end{align}
Moreover, for all $\delta\in (0,1)$, $k$ satisfies the following estimate:
\begin{align}\label{k_est1}
	\displaystyle|k(\zeta,t)|\leq |\zeta|\exp\left[\frac{1}{\min\{t,T-t\}}\left(\frac{\zeta^2}{\delta}-\frac{\beta}{1+\delta}\right)\right].
\end{align}
\end{lemma}

\noindent With the help of Lemma \ref{lemmaKernel}, we can prove the following observability result.

\begin{proposition}\label{prop62}
Let $s\in(0,1)$, $\partial\Omega^+$ as in \eqref{partition}, and fix $J\in\NN$. For any $T>0$, let $v_T\in \mathcal H_J$ and $v$ be the solution of the following adjoint backward parabolic system:
\begin{equation}\label{e61}
	\begin{cases}
		-v_t +(-\Delta)^sv=0 &\mbox{in }\;\Omega\times (0,T),
		\\
		v=0 &\mbox{in }\; (\mathbb R^N\setminus\Omega)\times (0,T),
		\\
		v(\cdot,T)=v_T &\mbox{in }\;\Omega.
	\end{cases}
\end{equation}
Then, there is a constant $C(s,\Omega)>0$ such that
\begin{align}\label{Ex1}
	\|v(\cdot,0)\|_{L^2(\Omega)}^2\le \frac{\Gamma(1+s)^2}{3sT}\exp\left(\frac{C(s,\Omega)}{T}T_0(J)^2\right)\int_0^T\int_{\partial\Omega^+}(x\cdot \nu)\left|\frac{v}{\rho^s}\right|^2\,d\sigma dt,
\end{align}
with $T_0(J)$ given by \eqref{e39}.
\end{proposition}

\noindent The proof of Proposition \ref{prop62} is given in Appendix \ref{appendix_proofs}.

\section{Controllability of the fractional heat equation} \label{sec7}

As a consequence of the observability inequality \eqref{Ex1}, we can now obtain our controllability results for the fractional heat equation. 

We start by showing the null controllability of the orthogonal projection of solutions to \eqref{heat_large} onto the space $\mathcal H_J$. Here solutions of \eqref{heat_large} are understood in the transposition sense introduced in Proposition 
\ref{prop:energy_parabolic}.

\begin{theorem}\label{controlThmHeat_bd}
Let $\partial\Omega^+$ be as in \eqref{partition} and $s\in(1/2,1)$. For all $T>0$, $u_0\in H^{-s}(\Omega)$, and $J\in\NN$ fixed, there exists a control function $f\in L^2(\partial\Omega^+\times(0,T))$ such that the solution $u$ of \eqref{heat_large} satisfies
\begin{align*}
	\Pi_{\mathcal H_J}u(\cdot,T)=0 \;\mbox{ a.e. in }\;\Omega.
\end{align*} 
Moreover, we have the following estimate of the controllability cost
\begin{align}\label{cost}
	\|f\|_{L^2(\partial\Omega^+\times(0,T))}^2\le \frac{C(s,\Omega)}{T}\exp\left(\frac{C(s,\Omega)}{T}T_0(J)^2\right)\|u_0\|_{H^{-s}(\Omega)}^2,
\end{align}
with $T_0(J)$ given by \eqref{e39}. If, furthermore, $u_0\in L^2(\Omega)$, the same controllability result holds with controllability cost 
\begin{align*}
	\|f\|_{L^2(\partial\Omega^+\times(0,T))}^2\le \frac{C(s,\Omega)}{T}\exp\left(\frac{C(s,\Omega)}{T}T_0(J)^2\right)\|u_0\|_{L^2(\Omega)}^2.
\end{align*}
\end{theorem}

\begin{remark}
\em{
While the result is expected to hold by duality, handling duality in the present fractional setting is delicate due to the non-standard nature of the boundary conditions in \eqref{heat_large}. This issue requires and will receive special attention in our analysis.
}
\end{remark}

\begin{proof}[\bf Proof of Theorem \ref{controlThmHeat_bd}]
		
We first treat the case of initial data in $H^{-s}(\Omega)$, from which the result for $L^2(\Omega)$ follows immediately.

\medskip
\noindent\textbf{Step 1: the case $u_0\in H^{-s}(\Omega)$.}
Since $\mathcal H_J\subset H_0^s(\Omega)$ is finite-dimensional, the condition
\begin{align*}
	\Pi_{\mathcal H_J}u(\cdot,T)=0 \quad\mbox{ in } H^{-s}(\Omega)
\end{align*}
is equivalent to
\begin{align*}
	\left\langle \Pi_{\mathcal H_J}u(\cdot,T),v_T\right\rangle_{-s,s}=0 \quad\text{for all } v_T\in\mathcal H_J.
\end{align*}

Let $v$ be the solution of \eqref{e61} associated with $v_T\in\mathcal H_J$. Multiplying \eqref{heat_large} by $v$ and integrating by parts using \eqref{eq:IBP_large}, we obtain
\begin{align*}
	\left\langle \Pi_{\mathcal H_J}u(\cdot,T),v_T\right\rangle_{-s,s}
	=
	\mathcal A(s)\int_0^T\int_{\partial\Omega^+} f\frac{v}{\rho^s}\,d\sigma dt+\left\langle u_0, v(\cdot,0)\right\rangle_{-s,s}.
\end{align*}
Consequently, the null controllability condition is equivalent to
\begin{align}\label{controlID_heat}
	\mathcal A(s)\int_0^T\int_{\partial\Omega^+} f\frac{v}{\rho^s}\,d\sigma dt+\left\langle u_0, v(\cdot,0)\right\rangle_{-s,s}=0.
\end{align}
We then consider the functional
\begin{align*}
	F(v_T):=\frac{\mathcal A(s)}{2}\int_0^T\int_{\partial\Omega^+}(x\cdot\nu)\left|\frac{v}{\rho^s}\right|^2\,d\sigma dt+\left\langle u_0,v(\cdot,0)\right\rangle_{-s,s}.
\end{align*}

The observability inequality \eqref{Ex1} implies that $F$ is convex and coercive over $\mathcal H_J$. Moreover, since $v(\cdot,0)\in\mathcal H_J$, we have
\begin{align*}
	\|v(\cdot,0)\|_{H_0^s(\Omega)}^2 \leq \lambda_J\|v(\cdot,0)\|_{L^2(\Omega)}^2.
\end{align*}
Hence,
\begin{align*}
	\left|\left\langle u_0,v(\cdot,0)\right\rangle_{-s,s}\right| \leq \|u_0\|_{H^{-s}(\Omega)}\|v(\cdot,0)\|_{H_0^s(\Omega)} \leq \lambda_J^{1/2}\|u_0\|_{H^{-s}(\Omega)}\|v(\cdot,0)\|_{L^2(\Omega)}.
\end{align*}

Therefore $F$ is also continuous on $\mathcal H_J$, so it admits a unique minimizer $v_T^\ast\in\mathcal H_J$. The corresponding Euler-Lagrange equation 
\begin{align}\label{EL_heat_Hs}
	\mathcal A(s)\int_0^T\int_{\partial\Omega^+}(x\cdot\nu)\frac{v}{\rho^s}\frac{v^\ast}{\rho^s}\,d\sigma dt + \langle u_0,v(\cdot,0)\rangle_{-s,s} = 0, \quad\text{ for all } v_T\in \mathcal H_J,
\end{align}
gives a control
\begin{align*}
	f = (x\cdot\nu)\frac{v^\ast}{\rho^s}\chi_{\partial\Omega^+}.
\end{align*}
Applying \eqref{EL_heat_Hs} to $v_T^\ast$ and using \eqref{Ex1}, we then get
\begin{align*}
	\bigg|\frac{1}{\mathcal A(s)} \langle u_0, v^\ast(\cdot,0)\rangle_{-s,s}\bigg| & \leq \frac{\lambda_J^{1/2}}{\mathcal A(s)}\|u_0\|_{H^{-s}(\Omega)}\|v(\cdot,0)\|_{L^2(\Omega)} \notag
	\\
	&\leq\frac{C(s)}{\sqrt{T}}\lambda_J^{1/2}\exp\left(\frac{C(s,\Omega)}{T}T_0(J)^2\right)\|u_0\|_{H^{-s}(\Omega)}\left(\int_0^T\int_{\partial\Omega^+}(x\cdot\nu)\left|\frac{v^\ast}{\rho^s}\right|^2\,d\sigma dt\right)^{\frac 12}.
\end{align*}

Moreover, since $T_0(J)=C\lambda_J^\gamma$ with $\gamma>0$, we have $T_0(J)^2\sim \lambda_J^{2\gamma}$. In particular, there exists a constant $C(s,\Omega)>0$ such that
\begin{align*}
	\lambda_J^{1/2}\leq \exp\left(C\lambda_J^{2\gamma}\right)\leq \exp\left(\frac{C}{T}T_0(J)^2\right),
\end{align*}
after enlarging the constant $C=C(s,\Omega)$ if necessary. Hence, the polynomial factor $\lambda_J^{1/2}$ can be absorbed into the exponential term, and we obtain
\begin{align}\label{cost_est}
	\bigg|\frac{1}{\mathcal A(s)} \langle u_0, v^\ast(\cdot,0)\rangle_{-s,s}\bigg| \leq\frac{C(s)}{\sqrt{T}}\exp\left(\frac{C(s,\Omega)}{T}T_0(J)^2\right)\|u_0\|_{H^{-s}(\Omega)}\left(\int_0^T\int_{\partial\Omega^+}(x\cdot\nu)\left|\frac{v^\ast}{\rho^s}\right|^2\,d\sigma dt\right)^{\frac 12}.
\end{align}
Combining \eqref{EL_heat_Hs} applied to $v^\ast$ and \eqref{cost_est} yields 
\begin{align*}
	\left(\int_0^T\int_{\partial\Omega^+}(x\cdot\nu)\left|\frac{v^\ast}{\rho^s}\right|^2\,d\sigma dt\right)^{\frac 12} \leq \frac{C(s)}{\sqrt{T}}\exp\left(\frac{C(s,\Omega)}{T}T_0(J)^2\right)\|u_0\|_{H^{-s}(\Omega)}.
\end{align*}
Therefore, 
\begin{align*}
	\|f\|_{L^2(\partial\Omega^+\times(0,T))}^2 &= \int_0^T\int_{\partial\Omega^+}(x\cdot\nu)^2\left|\frac{v^\ast}{\rho^s}\right|^2\,d\sigma dt \leq C(\Omega)\int_0^T\int_{\partial\Omega^+}(x\cdot\nu)\left|\frac{v^\ast}{\rho^s}\right|^2\,d\sigma dt 
	\\
	&\leq \frac{C(s,\Omega)}{T}\exp\left(\frac{C(s,\Omega)}{T}T_0(J)^2\right)\|u_0\|_{H^{-s}(\Omega)}^2.
\end{align*}

\smallskip
\noindent\textbf{Step 2: the case $u_0\in L^2(\Omega)$.}
This follows exactly as in the previous argument, replacing the duality pairing $\langle\cdot,\cdot\rangle_{-s,s}$ by the scalar product in $L^2(\Omega)$. In this case one obtains directly
\begin{align*}
	\|f\|_{L^2(\partial\Omega^+\times(0,T))}^2\le \frac{C(s,\Omega)}{T}\exp\left(\frac{C(s,\Omega)}{T}T_0(J)^2\right)\|u_0\|_{L^2(\Omega)}^2,
\end{align*}
and the proof is complete.
\end{proof}

\begin{remark}
\em{
Let us stress that, in \eqref{controlID_heat}, the constant $\mathcal A(s)$ is inherited from the integration by parts formula \eqref{eq:IBP_large} in which it appears as a normalization constant allowing us to recover the classical Green formula for the Laplacian operator at the limit, as $s\to 1^-$. Of course, this constant $\mathcal A(s)$ could be absorbed into the control, to get a clean duality identity
\begin{align*}
	\int_0^T\int_{\partial\Omega^+} f\frac{v}{\rho^s}\,d\sigma dt + \langle u_0, v(\cdot,0)\rangle_{-s,s} = 0,
\end{align*}
and functional 
\begin{align*}
	F(v_T)=\frac 12\int_0^T\int_{\partial\Omega^+}\left|\frac{v}{\rho^s}\right|^2\,d\sigma dt + \langle u_0, v(\cdot,0)\rangle_{-s,s}.
\end{align*} 
But for consistency with \eqref{eq:IBP_large}, we have decided to keep explicitly $\mathcal A(s)$ outside of the boundary integrals.
}
\end{remark}

Starting from the partial boundary controllability result established in Theorem \ref{controlThmHeat_bd}, we now prove our main result, Theorem \ref{ControlThmHeatIntro}. The proof relies on an iterative argument inspired by the approach in \cite{Le-Ro} for the classical heat equation (i.e., when $s = 1$), adapted to account for the specific features of the fractional setting considered here and the nature of low-frequency estimates and control results we achieved.

\begin{proof}[\bf Proof of Theorem \ref{ControlThmHeatIntro}]	
The starting point is to introduce a partition of the time interval $[0,T]$ defined as 
\begin{align*}
	[0,T] = \bigcup_{j\in\NN} [a_j,a_{j+1}],
\end{align*}
with 
\begin{align*}
	a_0:=0,\quad a_{j+1}:=a_j+2\tau_j, 
\end{align*}
and where
\begin{align*}
	\tau_j\coloneqq \gamma2^{-j\alpha},
\end{align*}
for fixed exponents 
\begin{align*}
	\varepsilon\in\left(0,s-\frac 12\right) \quad\text{ and }\quad 0<\alpha<\alpha^\ast\coloneqq\frac{2s-1-2\varepsilon}{2(1-2\varepsilon)}.
\end{align*}
The constant $\gamma>0$ is chosen so that
\begin{align}\label{LLL}
	2\sum_{j\in\NN} \tau_j = T=\lim_{j\to +\infty} a_j.
\end{align}
We now define the space 
\begin{align}\label{projection_space_heat}
	\mathcal K_j := \mbox{span}\Big\{\phi_1,\ldots,\phi_\ell\,:\,\lambda_\ell\leq \Lambda_j\Big\}, \quad\text{ with } \Lambda_j=2^j \text{ for all } j\in\NN,
\end{align}
and notice that $\mathcal K_j=\mathcal H_{J_j}$, where 
\begin{align*}
	J_j\coloneqq\max\{\ell\in\NN\,:\,\lambda_\ell\leq \Lambda_j\}.
\end{align*}
In particular, Theorem \ref{controlThmHeat_bd} applies to $\mathcal K_j$, and since $\lambda_{J_j}\leq \Lambda_j$, for $s>1/2$ we have
from \eqref{e39} that 
\begin{align*}
	T_0(J_j)\le C\Lambda_j^{\frac{1-s}{1-2\varepsilon}}.	
\end{align*}
We now construct the control function $f$ in the following way.
\begin{itemize}
	\item If $t\in(a_j,a_j+\tau_j]$, thanks to Theorem \ref{controlThmHeat_bd}, we apply a $L^2(\partial\Omega^+\times(0,T))$-control such that 
	\begin{align*}
		\Pi_{\mathcal K_j} u(\cdot,a_j+\tau_j)=0 \quad\text{ a.e. in } \Omega.	
	\end{align*}
	Since we know from Proposition \ref{prop:energy_parabolic} that $u\in C([0,T];H^{-s}(\Omega))$, we can use the energy estimate \eqref{eq:energy_parabolic_C0}, together with \eqref{cost}, to estimate
	\begin{align}\label{normEst1}
		\|u(\cdot,a_j+\tau_j)\|_{H^{-s}(\Omega)} &\leq C(s,\Omega)\Big[\|u(\cdot,a_j)\|_{H^{-s}(\Omega)} + \|f\|_{L^2(\partial\Omega^+\times(a_j,a_j+\tau_j))}\Big] \notag
		\\
		&\leq \frac{C(s,\Omega)}{\sqrt{\tau_j}}\exp\left(\frac{C(s,\Omega)}{\tau_j}T_0(J_j)^2\right)\|u(\cdot,a_j)\|_{H^{-s}(\Omega)} \notag
		\\
		&\leq \frac{C(s,\Omega,\varepsilon)}{\sqrt{\gamma}}2^{\frac{j\alpha}{2}}\exp\left(\frac{C(s,\Omega,\varepsilon)}{\gamma}2^{j\left(\frac{2(1-s)}{1-2\varepsilon}+\alpha\right)}\right)\|u(\cdot,a_j)\|_{H^{-s}(\Omega)}.
	\end{align}
	
	\item On the interval $(a_j + \tau_j, a_{j+1}]$, we let the system evolve freely, in the absence of control, to take advantage of the exponential decay of the solution. Thus,
	\begin{align}\label{normEst2}
		\|u(\cdot,a_{j+1})\|_{H^{-s}(\Omega)} & \leq \exp\left(-\tau_j2^j\right)\|u(\cdot,a_j+\tau_j)\|_{H^{-s}(\Omega)} \notag
		\\
		&= \exp\left(-\gamma 2^{j(1-\alpha)}\right)\|u(\cdot,a_j+\tau_j)\|_{H^{-s}(\Omega)}.
	\end{align}
\end{itemize}
Combining \eqref{normEst1} and \eqref{normEst2} we get
\begin{align*}
	\|u(\cdot,a_{j+1})\|_{H^{-s}(\Omega)} &\leq \frac{C(s,\Omega,\varepsilon)}{\sqrt{\gamma}}2^{\frac{j\alpha}{2}}\exp\left(\frac{C(s,\Omega,\varepsilon)}{\gamma} 2^{j\left(\frac{2(1-s)}{1-2\varepsilon}+\alpha\right)}-\gamma 2^{j(1-\alpha)}\right)\|u(\cdot,a_j)\|_{H^{-s}(\Omega)}.
\end{align*}
Iterating this procedure we obtain
\begin{align}\label{umb}
	\|u(\cdot,a_{j+1})\|_{H^{-s}(\Omega)} &\leq \frac{C(s,\Omega,\varepsilon)}{\sqrt \gamma}\exp\left(\sum_{\ell=0}^j \left(\frac{C(s,\Omega,\varepsilon)}{\gamma} 2^{\ell\left(\frac{2(1-s)}{1-2\varepsilon}+\alpha\right)}-\gamma 2^{\ell(1-\alpha)}\right)\right)\|u(\cdot,a_0)\|_{H^{-s}(\Omega)}\notag
	\\
	&\leq \frac{C(s,\Omega,\varepsilon)}{\sqrt \gamma}\exp\left(\sum_{\ell=0}^j \left(\frac{C(s,\Omega,\varepsilon)}{\gamma} 2^{\ell\left(\frac{2(1-s)}{1-2\varepsilon}+\alpha\right)}-\gamma 2^{\ell(1-\alpha)}\right)\right)\|u_0\|_{L^2(\Omega)}.
\end{align}
Notice that by \eqref{LLL}, since $u\in C([0,T];H^{-s}(\Omega))$, we have that
\begin{align*}
	\lim_{j\to +\infty} \|u(\cdot,a_{j+1})\|_{H^{-s}(\Omega)} = \|u(\cdot,T)\|_{H^{-s}(\Omega)}.	
\end{align*}
Moreover, for $\alpha<\alpha^\ast$, we also have that
\begin{align*}
	\frac{2(1-s)}{1-2\varepsilon}+\alpha < 1-\alpha.
\end{align*}

Hence, the sum in the exponential term in \eqref{umb} diverges to $-\infty$ as $j\to +\infty$ and we can deduce that $u(\cdot,T)=0$ a.e. in $\Omega$. 

Finally, the control $f$ built piecewise as above belongs to $L^2(\partial\Omega^+\times(0,T))$ and fulfills estimate \eqref{control_cost}. In fact, we can write 
\begin{align*}
	\|f \|_{L^2(\partial\Omega^+\times(0,T))}^2 = \sum_{j\in\NN} \|f\|_{L^2(\partial\Omega^+\times(a_j,a_j+\tau_j))}^2,
\end{align*}
where on each interval $(a_j,a_j+\tau_j)$, thanks to Theorem \ref{controlThmHeat_bd} we have 
\begin{align*}
		\|f\|_{L^2(\partial\Omega^+\times(a_j,a_j+\tau_j))}^2 &\leq \frac{C(s,\Omega)}{\tau_j}\exp\left(\frac{C(s,\Omega)}{\tau_j}T_0(J_j)^2 \right)\|u(\cdot,a_j)\|_{H^{-s}(\Omega)}^2
		\\
		&\leq \frac{C(s,\Omega,\varepsilon)}{\gamma} 2^{j\alpha} \exp\left(\frac{C(s,\Omega,\varepsilon)}{\gamma} 2^{j\left(\frac{2(1-s)}{1-2\varepsilon}+\alpha\right)}\right)\|u(\cdot,a_j)\|_{H^{-s}(\Omega)}^2.
\end{align*}
Combining this with the decay estimate
\begin{align*}
	\|u(\cdot,a_j)\|_{L^2(\Omega)}^2 \leq C\exp\left(-c\gamma 2^{j(1-\alpha)}\right)\|u_0\|_{L^2(\Omega)}^2,	
\end{align*}
and summing over $j$, we get
\begin{align*}
	\|f\|_{L^2(\partial\Omega^+\times(0,T))}^2 \leq \frac{C(s,\Omega,\varepsilon)}{\gamma} \sum_{j\in\NN} 2^{j\alpha} \exp\left(\frac{C(s,\Omega,\varepsilon)}{\gamma}2^{j\left(\frac{2(1-s)}{1-2\varepsilon}+\alpha\right)} - c\gamma 2^{j(1-\alpha)}\right)\|u_0\|_{L^2(\Omega)}^2.
\end{align*}

Since $0<\alpha<\alpha^\ast$ and $s\in(1/2,1)$, the above series converges and we get $f\in L^2(\partial\Omega^+\times(0,T))$. As for the control cost, notice that the argument in the exponential
\begin{align*}
	\frac{C(s,\Omega,\varepsilon)}{\gamma}2^{j\left(\frac{2(1-s)}{1-2\varepsilon}+\alpha\right)} - c\gamma 2^{j(1-\alpha)}
\end{align*}
reaches its maximum when 
\begin{align*}
	2^j\sim \gamma^{-\frac{2}{1-\alpha-\left(\frac{2(1-s)}{1-2\varepsilon}+\alpha\right)}} = \gamma^{-\frac{2(1-2\varepsilon)}{2s-1-2\varepsilon-2\alpha(1-2\varepsilon)}}.
\end{align*}
Substituting this into the positive term gives
\begin{align*}
	\frac{C(s,\Omega,\varepsilon)}\gamma 2^{j\left(\frac{2(1-s)}{1-2\varepsilon}+\alpha\right)} \sim \gamma^{-\theta}, \quad \text{ with }\theta = 
	\frac{3-2s-2\varepsilon}{2s-1-2\varepsilon-2\alpha(1-2\varepsilon)}.
\end{align*}
for all fixed $0<\alpha<\alpha^\ast$. Using this and the definitions of $\gamma$ and $\tau_j$, we finally get
\begin{align*}
	\|f\|_{L^2(\partial\Omega^+\times(0,T))}^2 \leq C(s,\Omega)\exp\left(\frac{C(s,\Omega)}{T^{\theta}}\right)\|u_0\|_{L^2(\Omega)}^2.	
\end{align*}
\end{proof}

\subsection{Optimality of the assumption $s\in (1/2,1)$}\label{sec:optimality}

Our main controllability result, namely Theorem \ref{ControlThmHeatIntro}, is established under the restriction $s\in (1/2,1)$ on the order of the fractional Laplacian. While this may initially seem like a technical condition tied to the specific techniques used in our proofs, we believe that this range is essentially sharp, and that controllability of solutions to \eqref{heat_large} should generally fail when $s\in(0,1/2]$. Several key reasons support this viewpoint, as we explain below.
\begin{itemize}
	\item For low-frequency components of \eqref{heat_large}, the controllability cost grows exponentially in the eigenvalue $\lambda_j$, with a rate governed by the exponent $\gamma$ given in \eqref{e39}. The iterative control strategy we applied relies on counteracting this exponential growth via the decay of the fractional heat semigroup. This balance is only achievable when $s \in (1/2,1)$. In contrast, for $s\in (0,1/2]$, we have $\gamma \geq 1$, leading to prohibitively high controllability costs that the semigroup decay cannot offset.
	
	\item This frequency-dependent behavior arises directly from our use of multiplier techniques, particularly the estimate of the cross terms in
	\begin{align*} 
		\int_\Omega p_t (x \cdot \nabla p)\,dx.
	\end{align*}
	This estimate relies on sharp regularity properties of eigenfunctions of the fractional Laplacian, tailored to each regime of $s$. Given the precision of these results, substantial improvements appear unlikely within the present approach.
	
	\item Our approach also uses the parabolic trace regularity result in Proposition \ref{prop:trace_parabolic}, available for solutions generated by $L^2$ interior source terms when $s\in (1/2,1)$. On the other hand, as shown in \cite{abels2023fractional,fernandez2024integro}, for $s \in (0,1/2]$ the boundary trace generally fails to exist even as a distribution.
\end{itemize}

Lastly, we note that the assumption $s\in (1/2,1)$ is recurrent in the literature on controllability of fractional heat-type equations-either involving the integral or spectral fractional Laplacian (see \cite{ABPWZ,biccari2017controllability,BiWaZu,micu2006controllability,miller2006controllability,warma2018null}). This further suggests that the condition is not merely technical, but rather reflects a fundamental threshold for the controllability properties of such systems.

\section{Concluding remarks and open problems}\label{sec-con-rem}

We have established boundary null controllability for the heat equation driven by the integral fractional Laplacian on a bounded smooth domain, for every $T>0$ and $s\in(1/2,1)$. The proof relies on a frequency-dependent boundary observability estimate for the associated fractional wave equation, transferred to the parabolic setting through transmutation and combined with the high-frequency dissipation of the fractional heat semigroup through a Lebeau-Robbiano iteration. A distinctive feature of the analysis is that the observability time for the fractional wave equation increases with the spectral cutoff, reflecting the deterioration of propagation at high frequencies.

Beyond the controllability result itself, our analysis also highlights several features of the fractional dynamics and its relation with the classical local case. In particular, we observe the following facts.
\begin{itemize}
	\item The controllability properties obtained for the fractional heat equation \eqref{heat_large} take place for all $T>0$, with a controllability cost which blows-up exponentially as $T\downarrow 0^+$. This is in accordance with the analogous properties of the local heat equation ($s=1$). See, for instance, \cite{fernandez2000cost,fernandez2000null} for more details.
	
	\item When passing to the limit formally as $s\to 1^-$ we get
	\begin{align*}
		f = \left.\frac{v^\ast_s}{\rho^s}\right|_{\partial\Omega^+}\overset{s\to 1}{\longrightarrow} \left.\frac{v^\ast}{\rho}\right|_{\partial\Omega^+} = \left.\frac{\partial v^\ast}{\partial\nu}\right|_{\partial\Omega^+},
	\end{align*}
	where we have used the notation $v^\ast_s$ to highlight the dependence on $s$ of the solutions to the adjoint equation \eqref{e61}. That is, our result is consistent with the classical literature for the local heat equation, for which we know that the boundary controls are given by the normal derivative of the adjoint state. Making this formal argument rigorous requires further analysis (see \cite{biccari2018poisson} for related results).  
	
	\item The duality argument used in the proof of Theorem \ref{controlThmHeat_bd} can similarly be applied starting from the observability inequality \eqref{e37} to establish the controllability of the projection of solutions onto $\mathcal{H}_J$ for the fractional wave equation. However, this does not yield information about the high-frequency components. As discussed in the previous sections, full null controllability in finite time cannot be expected, since the minimal controllability time $T_0$ tends to infinity as $J \to +\infty$.
\end{itemize}

Finally, several open problems and future research directions arise from our work, some of which are described below.

\begin{itemize}
	\item[1.] \textbf{Carleman estimates for the fractional Laplacian on bounded domains}. The transmutation approach used here does not apply to equations with space-time dependent potentials, limiting in particular its use for semilinear problems. An important alternative would be to develop Carleman estimates for the fractional heat equation on bounded domains. To the best of our knowledge, such estimates are currently available only in the whole space \(\mathbb{R}^N\) (see, e.g., \cite{ruland2015unique}), while their extension to bounded domains remains open due to the nonlocal character of the operator.
	
	\item[2.] \textbf{Optimality of the controllability region.} Our approach yields observability and controllability from the region $\partial\Omega_+$ where $x\cdot\nu>0$, which arises naturally from the multiplier argument for the fractional wave equation and is closely related to the Geometric Control Condition introduced in \cite{bardos1992sharp}. By contrast, the classical heat equation is controllable from any nonempty open subset of the boundary. Whether the same property holds for the fractional heat equation remains an open question beyond the techniques developed here.
	
	\item[3.] \textbf{Optimality of the observation time $T_0$.} Proposition \ref{prop31} gives $T_0(J)=C\lambda_J^{\gamma(s)}$, with the exponent $\gamma(s)$ in \eqref{gammarange} depending discontinuously on $s$. Since this discontinuity results from the different eigenfunction regularity estimates used in the proof, it is natural to ask whether sharper arguments could yield a continuous, or even optimal, dependence of $\gamma$ on $s$. We do not expect such an improvement to change the controllability threshold $s>1/2$ discussed in Section \ref{sec:optimality}.
	
	\item[4.] \textbf{Micro-local analysis for the fractional wave equation.} A natural direction is to investigate the frequency-dependent observability of the fractional wave equation using microlocal techniques. Such methods are fundamental for the classical wave equation (see, e.g., \cite{bardos1992sharp}), but, to the best of our knowledge, have not yet been developed in this context for the fractional Laplacian.	
\end{itemize}

\section{Acknowledgments}
The authors would like to express their sincere gratitude to Xavier Ros-Oton (University of Barcelona, Spain), Liviu Ignat (Polytechnic University of Bucharest and Simion Stoilow Institute of Mathematics of the Romanian Academy), and Nicola De Nitti (Polytechnic University of Bari, Italy) for valuable discussions and suggestions, which were instrumental in the proof of some of our results. 

{\appendix\section{Functional framework}\label{appendix}

We introduce here the functional framework and present some technical results on the fractional Laplacian which we used throughout the paper. 

\subsection{Formal definition and properties of the fractional Laplacian}
We start by giving a rigorous definition of the fractional Laplace operator. To this end, for any $s\in(0,1)$, we consider the space 
\begin{align*}
	\mathcal L^s(\RR^N):=\left\{u:\RR^N\to\mathbb R \;\mbox{ measurable such that }\;\int_{\Omega}\frac{|u(x)|}{(1+|x|)^{N+2s}}\;dx<+\infty\right\}.
\end{align*}
For $u\in \mathcal L^s(\mathbb R^N)$ and $\varepsilon>0$, we let
\begin{align*}
	(-\Delta)_\varepsilon^su(x):=C_{N,s}\int_{\{y\in\mathbb R^N:\;|x-y|\ge \varepsilon\}}\frac{u(x)-u(y)}{|x-y|^{N+2s}}\;dy,\;\;x\in\mathbb R^N,
\end{align*}
where the normalization constant is given by
\begin{align*}
	C_{N,s}:=\frac{s2^{2s}\Gamma\left(\frac{2s+N}{2}\right)}{\pi^{\frac N2}\Gamma(1-s)},
\end{align*}
and $\Gamma$ is the Euler Gamma function. The fractional Laplace operator $(-\Delta)^s$ is then defined for every $u\in \mathcal L^s(\mathbb R^N)$ by the formula
\begin{align}\label{FLO}
	(-\Delta)^su(x):=C_{N,s}\mbox{P.V.}\int_{\mathbb R^N}\frac{u(x)-u(y)}{|x-y|^{N+2s}}\;dy=\lim_{\varepsilon\downarrow 0}(-\Delta)_\varepsilon^su(x),\quad x\in\mathbb R^N,
\end{align}
provided that the limit exists for a.e. $x\in\mathbb R^N$.  We refer \cite{di2012hitchhiker} and their references for the class of functions for which the limit in \eqref{FLO} exists.

Let us now introduce the function spaces needed to investigate our problems, that is, the fractional order Sobolev spaces. In what follows, we will only provide the definitions and properties which are relevant for our results. More complete presentations can be found in several references, including but not limited to \cite{adams2003sobolev,di2012hitchhiker,Gris,lions1968problemes,War}.

Let $\Omega\subset\RR^N$ ($N\ge 1$) be an arbitrary open set. For any $s\in(0,1)$, we define the fractional order Sobolev space 
\begin{align}\label{Hs-def}
	H^s(\Omega):=\left\{u\in L^2(\Omega):\;\int_{\Omega}\int_{\Omega}\frac{|u(x)-u(y)|^2}{|x-y|^{N+2s}}\;dxdy< +\infty\right\}
\end{align}
and we endow it with the norm given by
\begin{align}\label{Hs-norm}
	\|u\|_{H^s(\Omega)}:=\left(\int_{\Omega}|u|^2\;dx+\int_{\Omega}\int_{\Omega}\frac{|u(x)-u(y)|^2}{|x-y|^{N+2s}}\;dxdy\right)^{\frac 12}.
\end{align}
Notice that $H^s(\mathbb R^N)$ can be defined by using the Fourier transform (see e.g. \cite{Gris}).
We let
\begin{align*}
	H_0^s(\Omega):=\Big\{u\in H^s(\mathbb R^N):\; u=0\;\mbox{ in }\mathbb R^N\setminus\Omega\Big\}=\Big\{u\in H^s(\mathbb R^N):\operatorname{supp}[u]\subset\overline{\Omega}\Big\}.
\end{align*}

We notice that if $0<s\ne 1/2<1$ and $\Omega$ is bounded and has a Lipschitz continuous boundary, then by \cite[Chapter 1]{Gris}, 
\begin{align*}
	H_0^s(\Omega)=\overline{\mathcal D(\Omega)}^{H^s(\mathbb R^N)},
\end{align*}
where $\mathcal D(\Omega)$ denotes the space of all continuous infinitely differentiable functions with compact support in $\Omega$.

We denote by $H^{-s}(\Omega) := (H_0^s(\Omega))^\ast$ the dual space of $H_0^s(\Omega)$ with respect to the pivot space $L^2(\Omega)$ so that the following continuous and compact embeddings hold: 
\begin{align*} 
	H_0^s(\Omega)\hookrightarrow L^2(\Omega)\hookrightarrow H^{-s}(\Omega).
\end{align*} 
If $s\ge 1$, then the above spaces are defined as in \cite{BWZ1,Gris} and the references therein.

Now we assume that $\Omega\subset\RR^N$ ($N\ge 1$) is a bounded open set. The following integration by parts formula is well-know. Let $u\in H_0^s(\Omega)$ be such that $(-\Delta)^su\in L^2(\Omega)$.  Then for every $v\in H_0^s(\Omega)$ we have that
\begin{align}\label{IPF}
	\int_{\Omega}v(-\Delta)^su\,dx=\frac{C_{N,s}}{2}\int_{\mathbb R^N}\int_{\mathbb R^N}\frac{(u(x)-u(y))(v(x)-v(y))}{|x-y|^{N+2s}}\,dxdy.
\end{align}
Using the notation
\begin{align}\label{nablas}
\int_{\Omega} |\nabla^s u|^2 \,  dx := \frac{C_{N,s}}{2} \int_{\mathbb{R}^N} \int_{\mathbb{R}^N} \frac{|u(x) - u(y)|^2}{|x - y|^{N + 2s}} \, dx \, dy,
\end{align}
we deduce from \eqref{IPF} that, for every $u \in H_0^s(\Omega)$ such that $(-\Delta)^s u \in L^2(\Omega)$, the following identity holds:
\begin{align*}
\int_{\Omega} |\nabla^s u|^2\, dx = \int_{\Omega} u\, (-\Delta)^s u\, dx.
\end{align*}

Finally, we introduce the realization in $L^2(\Omega)$ of the fractional Laplace operator with the zero Dirichlet exterior condition, that is, the operator $(-\Delta)_D^s$ on $L^2(\Omega)$ given by 
\begin{equation}\label{Ope}
	\begin{cases}
		D((-\Delta)_D^s)&:=\Big\{u\in H_0^s(\Omega):\; ((-\Delta)^su)|_{\Omega}\in L^2(\Omega)\Big\},
		\\
		(-\Delta)_D^su&=((-\Delta)^su)|_{\Omega}\;\mbox{ a.e. in }\;\Omega.
	\end{cases}
\end{equation}
A more rigorous definition of $(-\Delta)_D^s$ is given in \cite{Cl-Wa}.

The operator $(-\Delta)_D^s$ also defines an isomorphism from $H_0^s(\Omega)$ into $H^{-s}(\Omega)$. The following result is a direct consequence of the results proved  in \cite[Theorem 4.1]{G-JFA} (see  also \cite{JP-No}).

\begin{proposition}\label{grubb}
Let $\Omega\subset\RR^N$ ($N\ge 1$) be a bounded smooth open set and $s\in(0,1)$. Then, the following assertions hold.
\begin{itemize}
	\item[(a)] If $s\in(0,1/2)$, then $D((-\Delta)_D^s)= H_0^{2s}(\Omega)$.
	\item[(b)] If $s=1/2$, then $D((-\Delta)_D^s)\subset H^{1-\varepsilon}(\Omega)\cap H_0^s(\Omega)$ for all $\varepsilon\in(0,1)$.
	\item[(c)] If $s\in(1/2,1)$, then $D((-\Delta)_D^s)\subset H^{s+\frac 12-\varepsilon}(\Omega)\cap H_0^s(\Omega)$ for all $\varepsilon>0$.
\end{itemize}
\end{proposition}

We refer to \cite{Cl-Wa,G-JFA} and their references for more properties of $(-\Delta)_D^s$ and its  weak formulation. 

\subsection{PDEs involving the fractional Laplacian}

Let us now introduce some existence, uniqueness and regularity results for solutions of PDEs involving the fractional Laplacian, that we used  in our proofs.

\subsubsection{Fractional Poisson equation} Let us consider the following Dirichlet problem associated with the fractional Laplacian:
\begin{equation}\label{DP}
\begin{cases}
	(-\Delta)^su=f & \mbox{ in }\;\Omega,
	\\
	u=0 & \mbox{ in }\;\mathbb R^N\setminus\Omega.
\end{cases}
\end{equation}

The problem \eqref{DP} has been extensively studied in the last decade, and there are nowadays plenty of results about existence, uniqueness and regularity of solutions. In what follows, we present only a selection of such results, which are of relevance for our paper.

First of all, as a direct consequence of the classical Lax-Milgram Theorem, we have the following result about the existence and uniqueness of weak solutions to \eqref{DP} (see e.g. \cite[Proposition 2.1]{BWZ1} or \cite[Theorem 12]{leonori2015basic}).

\begin{proposition}\label{WPprop}
Let $\Omega\subset \RR^N$ be an arbitrary bounded open set and $s\in(0,1)$. Then for every $f\in H^{-s}(\Omega)$, the Dirichlet problem \eqref{DP} has a unique weak solution $u\in H_0^s(\Omega)$ fulfilling the identity
\begin{align*}
	\frac{C_{N,s}}{2}\int_{\mathbb R^N}\int_{\mathbb R^N}\frac{(u(x)-u(y))(v(x)-v(y))}{|x-y|^{N+2s}}\;dxdy=\langle f,v\rangle_{-s,s},
\end{align*}
for every $v\in H_0^s(\Omega)$. In addition, there is a constant $C=C(N,s,\Omega)>0$ such that
\begin{align*}
	\|u\|_{H_0^s(\Omega)} \leq C\|f\|_{H^{-s}(\Omega)}.
\end{align*}
\end{proposition}

Secondly, in \cite{RS1,ros2014extremal} (see also \cite{ros2016nonlocal,ros2016regularity,servadei2014weak} and the references therein), the authors investigated the Lebesgue and H\"older regularity of solutions to \eqref{DP}. In particular, the following results are available.

\begin{proposition}[{\cite[Proposition 1.4]{ros2014extremal}}]\label{prop:DP_Lebesgue}
Let $\Omega\subset\RR^N$ be a bounded $C^{1,1}$ domain, $s\in (0,1)$, $N > 2s$, $f\in C(\overline\Omega)$ and $u$ be the solution of \eqref{DP}.
\begin{itemize}
	\item[(a).] For each $1\leq r < N/(N-2s)$ there is a constant $C=C(N,s,r,\Omega)>0$ such that 
	\begin{align*}
		\|u\|_{L^r(\Omega)}\leq C\|f\|_{L^1(\Omega)}.
	\end{align*}
	\item[(b).] Let $1<p< N/(2s)$. Then there is a constant $C=C(N,s,p)>0$ such that 
	\begin{align*}
		\|u\|_{L^q(\Omega)}\leq C\|f\|_{L^p(\Omega)}, \quad\text{ with } q=\frac{Np}{N-2ps}.
	\end{align*}
	\item[(c).] Let $N/(2s)<p<+\infty$. Then, there is a constant $C=C(N,s,p,\Omega)>0$ such that
	\begin{align*}
		\|u\|_{C_0^\beta(\RR^N)}\leq C\|f\|_{L^p(\Omega)}, \quad\text{ with } \beta = \min\left\{s,2s-\frac Np\right\}.
	\end{align*}
\end{itemize}
\end{proposition}

\begin{proposition}[{\cite[Proposition 1.1]{RS1}}]\label{prop:DP_Holder}
Let $\Omega\subset\RR^N$ be a bounded $C^{1,1}$ domain. Given $f\in L^\infty(\Omega)$, let $u$ be the solution of \eqref{DP}. Then, $u\in C_0^s(\RR^N)$ and there exists a constant $C=C(s,\Omega)>0$ such that $$\|u\|_{C_0^s(\RR^N)}\leq C\|f\|_{L^\infty(\Omega)}.$$
\end{proposition}

\begin{proposition}[{\cite[Theorem 1.2]{RS1}}]\label{prop:DP_Holder_boundary}
Let $\Omega\subset\RR^N$ be a bounded $C^{1,1}$ domain, $f\in L^\infty(\Omega)$, $u$ be the solution of \eqref{DP}, and $\rho(x) = \text{dist}(x,\partial\Omega)$. Then, $(u/\rho^s)|_\Omega$ can be continuously extended to $\overline\Omega$. Moreover, we have that $(u/\rho^s)\in C^\alpha(\overline\Omega)$ and 
\begin{align*}
	\left\|\frac{u}{\rho^s}\right\|_{C^\alpha(\overline\Omega)} \leq C(s,\Omega)\|f\|_{L^\infty(\Omega)}	
\end{align*}
for some $0<\alpha< \min\{s, 1-s\}$. 
\end{proposition}

In addition to that, we also have the following result for the trace of $u/\rho^s$ over $\partial\Omega$ (see \cite[Proposition 4.4]{abels2023fractional} and \cite[Section 2.8.3]{fernandez2024integro}).

\begin{proposition}\label{prop:trace}
Let $\Omega\subset\RR^N$ be a bounded $C^2$ domain. Given $f\in L^2(\Omega)$, let $u$ be the solution of \eqref{DP}. Then, for all $s\in(1/2,1)$, we have that $(u/\rho^s)|_{\partial\Omega}\in H^{s-1/2}(\partial\Omega)$ and 
\begin{align}\label{eq:trace_est}
	\left\|\frac{u}{\rho^s}\right\|_{H^{s-1/2}(\partial\Omega)} \leq C(\Omega)\|f\|_{L^2(\Omega)}.
\end{align}
\end{proposition}

All the previously mentioned results provide optimal H\"older regularity up to the boundary for weak solutions of the non-local Poisson problem \eqref{DP}. However, explicit computations reveal that there also exist pointwise solutions of \eqref{DP}, which explode at the boundary $\partial\Omega$, behaving asymptotically like $\rho^{s-1}$. These are the so-called \textit{large solutions} of \eqref{DP}, which have been studied, for instance, in \cite{abatangelo2015large,abatangelo2017very,abatangelo2023singular,chan2021blow}. In particular, we have the following result. 

\begin{proposition}[{\cite[Theorem 1.2.3]{abatangelo2015large}}]\label{prop:large_solutions}
Let $\Omega\subset \RR^N$ be a bounded $C^{1,1}$ domain. Given $\alpha\in (0,1)$, let $f\in C^\alpha(\overline\Omega)$, and $h\in C(\partial\Omega)$. Then, for any $s\in(0,1)$ and $\varepsilon >0$, there exists a unique function $u\in C^{2s+\varepsilon}_{\text{loc}}(\Omega)$ which is a point-wise solution of the Dirichlet problem 
\begin{equation}\label{DP_large}
	\begin{cases}
		(-\Delta)^su=f & \mbox{ in }\;\Omega,
		\\
		u=0 & \mbox{ in }\;\RR^N\setminus\overline\Omega,
		\\
		\displaystyle\rho^{1-s}u=h & \mbox{ on }\;\partial\Omega.
	\end{cases}
\end{equation}
\end{proposition}

Finally, we have the following integration by parts formula involving large and standard solutions of \eqref{DP} (see \cite[Proposition 1.2.2]{abatangelo2015large} or \cite[Lemma 5.3]{fernandez2024stable}).

\begin{proposition}\label{prop:IBP}
Let $\Omega\subset\RR^N$ be a bounded $C^2$ domain and $s\in(0,1)$. Let $u$ be a standard solution to \eqref{DP} and $v$ be a large solution fulfilling \eqref{DP_large}. Then,
\begin{align}\label{eq:IBP_large}
	\int_\Omega v(-\Delta)^su\,dx = \int_\Omega u(-\Delta)^sv\,dx + \Gamma(s)\Gamma(1+s)\int_{\partial\Omega} \frac{u}{\rho^s}\frac{v}{\rho^{s-1}}\,d\sigma. 
\end{align}
\end{proposition}

\subsubsection{Fractional heat equation} Let us consider the following heat-like equation associated with the fractional Laplacian:
\begin{equation}\label{heat_appendix}
	\begin{cases}
		u_t + (-\Delta)^s u = f & \mbox{in }\;\Omega\times (0,T),
		\\
		u=0 & \mbox{in }\; (\mathbb R^N\setminus\Omega)\times (0,T),
		\\
		u(\cdot,0)=u_0, &\mbox{in }\;\Omega.
	\end{cases}
\end{equation}

The existence and uniqueness of weak solutions given in Proposition \ref{WPprop} can be naturally extended also to this parabolic context. In particular, we have the following result.

\begin{proposition}[{\cite[Chapter 10, Theorem 10.9]{brezis2010functional}}]\label{weakSolHeat}
Let $\Omega\subset \RR^N$ be an arbitrary bounded open set and $s\in(0,1)$. For every $u_0\in L^2(\Omega)$ and $f\in L^2((0,T);H^{-s}(\Omega))$, the fractional heat equation \eqref{heat_appendix} admits a unique finite energy solution 
\begin{align*}
	u\in C([0,T];L^2(\Omega))\cap L^2((0,T);H_0^s(\Omega))\cap H^1((0,T);H^{-s}(\Omega))
\end{align*} 
such that $u(\cdot,0)=u_0$ a.e. in $\Omega$ and the identity
\begin{align}\label{fes}
	\int_0^T \langle f,w\rangle_{-s,s}\,dt = &\, \int_0^T \langle u_t,w\rangle_{-s,s}\,dt 
	\\
	&+ \frac{C_{N,s}}{2}\int_0^T\int_{\RR^N}\int_{\RR^N}\frac{(u(x,t)-u(z,t))(w(x)-w(z))}{|x-z|^{N+2s}}\,dxdzdt\notag
\end{align}
holds, for every $w\in H_0^s(\Omega)$.
\end{proposition}

\noindent We also have the result of maximal $L^2$-regularity taken from \cite[Theorem 1]{lamberton1987equations}.

\begin{proposition}\label{prop:maximal}
Let $\Omega\subset \RR^N$ be an arbitrary bounded open set and $s\in(0,1)$. Take $u_0=0$ and $f\in L^2(\Omega\times(0,T))$. Then, the fractional heat equation \eqref{heat_appendix} admits a solution 
\begin{align*}
	u\in C([0,T];L^2(\Omega))
\end{align*} 
such that 
\begin{align*}
	u_t, (-\Delta)^s u\in L^2(\Omega\times(0,T))
\end{align*}
and there is a constant $C>0$ (independent of $f$ and $u$) such that
\begin{align}\label{eq:maximal_estimate}
	\|u_t\|_{L^2(\Omega\times(0,T))} + \|(-\Delta)^su\|_{L^2(\Omega\times(0,T))} \leq C \|f\|_{L^2(\Omega\times(0,T))}.
\end{align}
\end{proposition}

Proposition \ref{prop:maximal}, combined with Proposition \ref{prop:trace}, allows us to establish a similar trace result for the fractional heat equation. 

\begin{proposition}\label{prop:trace_parabolic}
Let $\Omega\subset\RR^N$ be a bounded $C^2$ domain. Take $u_0=0$ and $f\in L^2(\Omega\times(0,T))$. Then, for all $s\in (1/2,1)$, the solution $u$ of \eqref{heat_appendix} satisfies 
\begin{align*}
	\left.\frac{u}{\rho^s}\right|_{\partial\Omega}\in L^2((0,T);H^{s-1/2}(\partial\Omega)), 	
\end{align*}
and there is a constant $C(\Omega)>0$ such that
\begin{align}\label{eq:trace_est_parabolic}
	\left\|\frac{u}{\rho^s}\right\|_{L^2((0,T);H^{s-1/2}(\partial\Omega))} \leq C(\Omega)\|f\|_{L^2(\Omega\times(0,T))}.
\end{align}
\end{proposition}

\begin{proof}
For all $t\in (0,T)$, we have from \eqref{heat_appendix} that $u(t)$ solves the elliptic problem
\begin{align*}
	\begin{cases}
		(-\Delta)^s u(t) = \psi(t) & \mbox{in }\;\Omega,
		\\
		u(t)=0 & \mbox{in }\; \mathbb R^N\setminus\Omega,
	\end{cases}
\end{align*}
with $\psi(t)\coloneqq f(t) - u_t(t)$. Moreover, because $f\in L^2(\Omega\times(0,T))$, we have from Proposition \ref{prop:maximal} that $u_t\in L^2(\Omega\times(0,T))$ and, therefore, $\psi\in L^2(\Omega\times(0,T))$. Since $s\in (1/2,1)$, using Proposition  \ref{prop:trace} and \eqref{eq:maximal_estimate}, we can then conclude that 
\begin{align*}
	\left\|\frac{u}{\rho^s}\right\|_{L^2((0,T);H^{s-1/2}(\partial\Omega))}^2 &= \int_0^T\left\|\frac{u(t)}{\rho^s}\right\|_{H^{s-1/2}(\partial\Omega)}^2\,dt \leq C(\Omega) \int_0^T \|\psi(t)\|_{L^2(\Omega)}^2\,dt 
	\\
	&\leq C(\Omega) \int_0^T \Big(\|f(t)\|_{L^2(\Omega)}+\|u_t(t)\|_{L^2(\Omega)}\Big)^2\,dt 
	\\
	&\leq C(\Omega) \int_0^T \|f(t)\|_{L^2(\Omega)}^2\,dt = C(\Omega)\|f\|_{L^2(\Omega\times(0,T))}^2,
\end{align*}
and the proof is finished.
\end{proof}

In \cite{fernandez2016boundary}, the authors established H\"older regularity up to the boundary for solutions to \eqref{heat_appendix}, in a form analogous to the results stated in Propositions \ref{prop:DP_Holder} and \ref{prop:DP_Holder_boundary}.

\begin{proposition}[{\cite[Theorem 1.1]{fernandez2016boundary}}]\label{prop:FHE_Holder}
Let $\Omega\subset\RR^N$ be a bounded $C^{1,1}$ domain and $s\in(0,1)$.  For $u_0\in L^2(\Omega)$, let $u$ be the unique weak solution of \eqref{heat_appendix} with $f\equiv 0$. Then, for all $T>0$ and $\varepsilon\in (0,s)$, there are two constants $C_1=C_1(N,s,\Omega,T)>0$ and $C_2=C_2(N,s,\Omega,T,\varepsilon)>0$ such that 
\begin{align*}
	&\sup_{t>0}\|u(\cdot,t)\|_{C^s(\RR^N)}\leq C_1\|u_0\|_{L^2(\Omega)},
	\\
	&\sup_{t>0}\left\|\frac{u(\cdot,t)}{\delta^s}\right\|_{C^{s-\varepsilon}(\overline\Omega)}\leq C_2\|u_0\|_{L^2(\Omega)}.
\end{align*}
Moreover, the constants $C_1$ and $C_2$ blow up as $T\downarrow 0^+$.	
\end{proposition}

Finally, analogously to what we have done for the elliptic equation \eqref{DP_large}, we can extend the notion of large solutions also to the fractional heat equation, by considering the parabolic model
\begin{equation}\label{heat_appendix_large}
	\begin{cases}
		u_t + (-\Delta)^s u = 0 & \mbox{in }\;\Omega\times (0,T),
		\\
		u=0 & \mbox{in }\; (\mathbb R^N\setminus\overline\Omega)\times (0,T),
		\\
		\rho^{1-s}u=h & \mbox{on }\; \partial\Omega\times (0,T),
		\\
		u(\cdot,0)=u_0, &\mbox{in }\;\Omega.
	\end{cases}
\end{equation}

Solutions to \eqref{heat_appendix_large} can be defined by transposition. In particular, in \cite[Theorems 1.1 and 8.3]{chan2022singular} it is shown that for all $s\in(0,1)$, initial datum $u_0\in L^1(\Omega,\rho^s dx)$, and boundary datum $f\in L^1(\partial\Omega\times(0,T))$, \eqref{heat_appendix_large} admits a unique weak-dual solution $u\in L^1(\Omega\times(0,T),\rho^s dxdt)$ satisfying the corresponding transposition identity against solutions of the backward homogeneous adjoint problem. 

When $s\in(1/2,1)$, thanks to Proposition \ref{prop:trace_parabolic}, we can get some improved regularity. To this end, it is convenient to separate the contribution of the initial datum from the contribution of the singular boundary datum. Accordingly, we write
\begin{align}\label{eq:u_decomp}
	u = \eta + \xi, \qquad \eta(t):=e^{-t(-\Delta)_D^s}u_0,
\end{align}
where $\xi$ solves
\begin{equation}\label{heat_appendix_large_zero}
	\begin{cases}
		\xi_t + (-\Delta)^s \xi = 0 & \mbox{in }\;\Omega\times (0,T),
		\\
		\xi=0 & \mbox{in }\; (\mathbb R^N\setminus\overline\Omega)\times (0,T),
		\\
		\rho^{1-s}\xi=h & \mbox{on }\; \partial\Omega\times (0,T),
		\\
		\xi(\cdot,0)=0, & \mbox{in }\;\Omega.
	\end{cases}
\end{equation}
We then have the following definition of transposition solution for \eqref{heat_appendix_large_zero}
\begin{definition}
Let $h\in L^2(\partial\Omega\times(0,T))$. We say that $\xi\in L^2(\Omega\times(0,T))$ is a transposition solution of \eqref{heat_appendix_large_zero} if, for every $q\in L^2(\Omega\times(0,T))$ and solution $z$ of the backward adjoint problem
\begin{equation}\label{heat_adjoint_large_zero} 
	\begin{cases}
		-z_t + (-\Delta)^s z = q & \mbox{in }\;\Omega\times (0,T),
		\\
		z=0 & \mbox{in }\; (\mathbb R^N\setminus\overline\Omega)\times (0,T),
		\\
		z(\cdot,T)=0, & \mbox{in }\;\Omega,
	\end{cases}
\end{equation}
it holds the identity 
\begin{align}\label{heat_large_transp}
	\int_0^T\int_\Omega \xi q\,dxdt = \Gamma(s)\Gamma(1+s)\int_0^T\int_{\partial\Omega} h\frac{z_q}{\rho^s}\,d\sigma dt.
\end{align}
\end{definition}

\begin{proposition}\label{prop:energy_parabolic}
Let $\Omega\subset\RR^N$ be a bounded $C^2$ domain, let $s\in (1/2,1)$, and let $h\in L^2(\partial\Omega\times(0,T))$. Then the following properties hold.

\begin{itemize}
	\item[(a)] For every $u_0\in H^{-s}(\Omega)$, \eqref{heat_appendix_large} admits a unique solution
	\begin{align*}
		u\in L^2(\Omega\times (0,T))\cap C([0,T];H^{-s}(\Omega)).
	\end{align*}
	Moreover, there exists a constant $C=C(s,\Omega,T)>0$ such that
	\begin{align}\label{eq:energy_parabolic_L2}
		\|u\|_{L^2(\Omega\times (0,T))}\leq C\Big(\|u_0\|_{H^{-s}(\Omega)} + \|h\|_{L^2(\partial\Omega\times (0,T))}\Big)
	\end{align}
	and
	\begin{align}\label{eq:energy_parabolic_C0}
		\|u\|_{C([0,T];H^{-s}(\Omega))}\leq C\Big(\|u_0\|_{H^{-s}(\Omega)} + \|h\|_{L^2(\partial\Omega\times (0,T))}\Big).
	\end{align}
	
	\item[(b)] If, in addition, $u_0\in L^2(\Omega)$, then the same transposition solution satisfies
	\begin{align*}
		\|u\|_{L^2(\Omega\times (0,T))}\leq C\Big(\|u_0\|_{L^2(\Omega)} + \|h\|_{L^2(\partial\Omega\times (0,T))}\Big)
	\end{align*}
	and
	\begin{align*}
		\|u\|_{C([0,T];H^{-s}(\Omega))}\leq C\Big(\|u_0\|_{L^2(\Omega)} + \|h\|_{L^2(\partial\Omega\times (0,T))}\Big).
	\end{align*}	
\end{itemize}
\end{proposition}

\begin{proof}
We divide the argument into three steps.

\medskip
\noindent\textbf{Step 1: existence, uniqueness and $L^2(\Omega\times(0,T))$ estimate for $\xi$.} 

Let $q\in L^2(\Omega\times(0,T))$, and let $z=z$ denote the solution of the backward adjoint problem \eqref{heat_adjoint_large_zero}. By the change of variables $\widetilde z(x,t):=z(x,T-t)$, the function $\widetilde z$ solves a forward fractional heat equation with zero initial datum and right-hand side $\widetilde q(x,t):=q(x,T-t)$. Hence, by Proposition \ref{prop:trace_parabolic},
\begin{align}\label{eq:adjoint_trace_est}
	\left\|\frac{z}{\rho^s}\right\|_{L^2((0,T);H^{s-1/2}(\partial\Omega))}\leq C(s,\Omega,T)\|q\|_{L^2(\Omega\times(0,T))}.
\end{align}
Therefore, the linear functional
\begin{align*}
	\mathcal L(q):=\Gamma(s)\Gamma(1+s)\int_0^T\int_{\partial\Omega} h\frac{z}{\rho^s}\,d\sigma dt
\end{align*}
is continuous on $L^2(\Omega\times(0,T))$, since
\begin{align*}
	|\mathcal L(q)| & \leq C(s)\|h\|_{L^2(\partial\Omega\times (0,T))} \left\|\frac{z}{\rho^s}\right\|_{L^2(\partial\Omega\times(0,T))}
	\\
	&\leq C(s,\Omega,T)\|h\|_{L^2(\partial\Omega\times (0,T))} \left\|\frac{z}{\rho^s}\right\|_{L^2((0,T);H^{s-1/2}(\partial\Omega))} \leq C(s,\Omega,T)\|h\|_{L^2(\partial\Omega\times (0,T))}\|q\|_{L^2(\Omega\times(0,T))}.
 \end{align*}
By the Riesz representation theorem, there exists a unique $\xi\in L^2(\Omega\times(0,T))$ satisfying \eqref{heat_large_transp}. Moreover,
\begin{align}\label{est_xi}
	\|\xi\|_{L^2(\Omega\times(0,T))} \leq C(s,\Omega,T)\|h\|_{L^2(\partial\Omega\times (0,T))}.
\end{align}

\medskip
\noindent\textbf{Step 2: estimates for the semigroup part $\eta$.}

Let $\eta(t)=e^{-t(-\Delta)_D^s}u_0$. If $u_0\in L^2(\Omega)$, then the contractivity of the semigroup in $L^2(\Omega)$ gives
\begin{align}\label{eq:eta_L2_est}
	\|\eta\|_{L^2(\Omega\times(0,T))}\leq \sqrt{T}\,\|u_0\|_{L^2(\Omega)}
	\quad\mbox{and}\quad
	\|\eta\|_{C([0,T];H^{-s}(\Omega))}\leq C(s,\Omega)\|u_0\|_{L^2(\Omega)}.
\end{align}
If $u_0\in H^{-s}(\Omega)$, then, writing
\begin{align*}
	u_0=\sum_{j\geq 1} u_{0,j}\phi_j, \qquad \eta(t)=\sum_{j\geq 1} e^{-\lambda_j t}u_{0,j}\phi_j,
\end{align*}
we obtain
\begin{align*}
	\|\eta\|_{L^2(\Omega\times(0,T))}^2
	= \sum_{j\geq 1}|u_{0,j}|^2\int_0^T e^{-2\lambda_j t}\,dt
	\leq \frac12\sum_{j\geq 1}\lambda_j^{-1}|u_{0,j}|^2
	\leq C(s,\Omega)\|u_0\|_{H^{-s}(\Omega)}^2,
\end{align*}
while
\begin{align}\label{eq:eta_Hms_est}
	\|\eta\|_{C([0,T];H^{-s}(\Omega))}\leq C(s,\Omega)\|u_0\|_{H^{-s}(\Omega)}.
\end{align}

\medskip
\noindent\textbf{Step 3: estimates and continuity in time for $u=\eta+\xi$.}
By construction, $\xi$ solves
\begin{align*}
	\xi_t = -(-\Delta)_D^s\xi
\end{align*}
in the sense of distributions. Since $(-\Delta)_D^s:L^2(\Omega)\to H^{-2s}(\Omega)$ is continuous, we get
\begin{align*}
	\|\xi_t\|_{L^2(0,T;H^{-2s}(\Omega))}\leq C(s,\Omega)\|\xi\|_{L^2(\Omega\times(0,T))}.
\end{align*}
Hence
\begin{align*}
	\xi\in L^2(\Omega\times(0,T))\cap H^1(0,T;H^{-2s}(\Omega)),
\end{align*}
and by interpolation (see \cite[Chapter 1, Section 3, Theorem 3.1]{lions1968problemes}) we conclude that
\begin{align*}
	\xi\in C([0,T];[L^2(\Omega),H^{-2s}(\Omega)]_{1/2})=C([0,T];H^{-s}(\Omega)),
\end{align*}
with
\begin{align}\label{eq:xi_C0_est}
	\|\xi\|_{C([0,T];H^{-s}(\Omega))}
	\leq C(s,\Omega,T)\Big(\|\xi\|_{L^2(\Omega\times(0,T))}+\|\xi_t\|_{L^2(0,T;H^{-2s}(\Omega))}\Big)
	\leq C(s,\Omega,T)\|h\|_{L^2(\partial\Omega\times(0,T))}.
\end{align}

If $u_0\in H^{-s}(\Omega)$, then \eqref{eq:energy_parabolic_L2} and \eqref{eq:energy_parabolic_C0} follow by combining \eqref{eq:u_decomp}, \eqref{est_xi}, the estimate above for $\eta$ in $L^2(\Omega\times(0,T))$, \eqref{eq:eta_Hms_est}, and \eqref{eq:xi_C0_est}. If $u_0\in L^2(\Omega)$, the corresponding bounds follow from \eqref{eq:u_decomp}, \eqref{est_xi}, \eqref{eq:eta_L2_est}, and \eqref{eq:xi_C0_est}. The uniqueness follows from the uniqueness of the transposition solution $\xi$.
\end{proof}

\begin{remark}
\em{
Proposition \ref{prop:energy_parabolic} ensures that the controlled fractional heat equation admits a unique transposition solution $u\in C([0,T];H^{-s}(\Omega))$. Consequently, the pointwise values $u(t)$ are well defined in $H^{-s}(\Omega)$ for every $t\in[0,T]$. In particular, the evaluations $u(a_j)$, $u(a_j+\tau_j)$, and $u(T)$ employed in the Lebeau-Robbiano iteration in the proof of Theorem \ref{ControlThmHeatIntro} are meaningful.
}
\end{remark}

\begin{remark}
\em{
The results we have presented in this section are stated under different regularity assumptions for the domain $\Omega$. Specifically, while some of these results like Propositions \ref{prop:DP_Lebesgue} or \ref{prop:DP_Holder} require $\Omega$ to be $C^{1,1}$, others like Propositions \ref{prop:trace} or \ref{prop:IBP} assume $C^2$-regularity, and some refer more vaguely to \textit{regular} domains without providing a precise definition (see, for instance, Proposition \ref{prop:maximal}). Throughout this work, in order to avoid repeating assumptions and to maintain consistency, we adopt $C^2$-regularity as a standing assumption. This level of regularity suffices for all our purposes, as it ensures the applicability of all results from the existing literature that are invoked in the paper. We emphasize that this choice is not intended to be optimal, but rather to streamline the exposition while encompassing all necessary cases.
}
\end{remark}

\subsection{Eigenvalues problem for the fractional Laplacian}

It is well known that the operator $(-\Delta)_D^s$, defined in \eqref{Ope}, has a compact resolvent. As a consequence, its spectrum is discrete and consists of a sequence of eigenvalues $(\lambda_j)_{j \in \mathbb{N}}$ satisfying
\begin{align*}
	0<\lambda_1\le\lambda_2\le\cdots\le\lambda_j\le\cdots\;\to +\infty, \mbox{ as } j\to+\infty.
\end{align*}

We denote by $(\phi_j)_{j\in\NN}$ the normalized eigenfunctions associated with the eigenvalues $(\lambda_j)_{j\in\NN}$, i.e., the solutions of the following Dirichlet problem:
\begin{equation}\label{e27}
	\begin{cases}
		(-\Delta)^s\phi_j=\lambda_j\phi_j & \mbox{in }\;\Omega,
		\\
		\phi_j=0 & \mbox{in }\; \mathbb R^N\setminus\Omega.
	\end{cases}
\end{equation}
For each eigenfunction $\phi_j$, normalized as
\begin{align}\label{e33}
	\|\phi_j\|_{L^2(\Omega)}^2 = \int_{\Omega}\phi_j^2\,dx=1,
\end{align}
we have that
\begin{align}\label{e34}
	\|\phi_j\|_{H_0^s(\Omega)}^2= \int_{\Omega}|\nabla^s\phi_j|^2\,dx=\lambda_j.
\end{align}

Moreover, we have the following Weyl law for the Dirichlet eigenvalues (see \cite[Proposition 2.1]{fernandez2016boundary} and also \cite{blumenthal1959asymptotic,FRGL}):
\begin{align*}
	\lambda_j\sim C(N,s,\Omega)j^{\frac{2s}{N}},
\end{align*}
in the sense that there are two constants $0<C_1(N,s,\Omega)\leq C_2(N,s,\Omega)$ such that
\begin{align}\label{eq:weyl}
	C_1(N,s,\Omega)j^{\frac{2s}{N}}\leq \lambda_j\leq C_1(N,s,\Omega)j^{\frac{2s}{N}}, \quad\text{ for all }j\in\NN.
\end{align}

Let us notice that, through the spectrum of $(-\Delta)_D^s$, it is also possible to characterize the operator's domain and its dual. More precisely, we have the following:
\begin{equation}\label{domain_char}
	\begin{array}{ll}
		\displaystyle D((-\Delta)^s_D) = \Big\{u\in L^2(\Omega):\displaystyle \|u\|_{D((-\Delta)^s_D)}^2 <+\infty\Big\}, & 
		\displaystyle \|u\|_{D((-\Delta)^s_D)}^2 = \sum_{j\geq 1} \left|\lambda_j(u,\phi_j)_{L^2(\Omega)}\right|^2,
		\\[20pt]
		\displaystyle D\big(((-\Delta)^s_D)^\ast\big)= \Big\{u\in L^2(\Omega):\displaystyle \|u\|_{D((-\Delta)^s_D)^\ast}^2<+\infty\Big\}, & \displaystyle \|u\|_{D((-\Delta)^s_D)^\ast}^2 = \sum_{j\geq 1} \left|\lambda_j^{-1}(u,\phi_j)_{L^2)\Omega)}\right|^2, 
	\end{array} 
\end{equation}
where $(\cdot,\cdot)_{L^2(\Omega)}$ denotes the scalar product in $L^2(\Omega)$.

In addition to that, when $s\in(1/2,1)$, Proposition \ref{grubb}(c) ensures that the eigenfunctions $(\phi_j)_{j\in\NN}$ belong to $H_0^1(\Omega)$. 

Finally, we give a series of important results on the eigenfunctions of the fractional Laplacian, which are used in the proof of the observability properties of the adjoint wave system \eqref{adjointWave}. 

\begin{proposition}\label{lem31}
Let $\Omega\subset\RR^N$ ($N\ge 1$) be a bounded smooth open set, let $s\in(1/2,1)$ and let $\varepsilon\in(0,s-\frac 12)$. Then, there is a constant $C=C(N,s,\Omega,\varepsilon)>0$ such that for every $j\in\NN$ the eigenfunction $\phi_j$ satisfies
\begin{align}\label{e35}
	\int_{\Omega}|\nabla \phi_j|^2\,dx\le C\lambda_j^{\frac{2(1-s)}{1-2\varepsilon}}\int_{\Omega}|\nabla^s\phi_j|^2\,dx.
\end{align}
\end{proposition}

\begin{proof}
Let us first notice that, since $s\in(1/2,1)$, we have from Proposition \ref{grubb}(c) that
\begin{align*}
	\phi_j\in D((-\Delta)_D^s)\subset H^{s+1/2-\varepsilon}(\Omega)\cap H_0^s(\Omega)\subset H_0^1(\Omega).
\end{align*}
Now, define the parameter
\begin{align*}
	\theta\coloneqq 1-\frac{2(1-s)}{1-2\varepsilon}\in (0,1),
\end{align*}
for which we have 
\begin{align*}
	1 = (1-\theta)\left(s+\frac 12-\varepsilon\right) + \theta s.
\end{align*}

Using the interpolation inequality between the Hilbert spaces $H^{s+1/2-\varepsilon}(\Omega)$, $H_0^1(\Omega)$ and $H_0^s(\Omega)$ (see \cite[Theorem 11.6]{lions1968problemes}), together with the fact that $D((-\Delta)_D^s)\hookrightarrow H^{s+1/2-\varepsilon}(\Omega)$, we get that there is a constant $C=C(N,s,\Omega,\varepsilon)>0$ such that for every $j\in\NN$ 
\begin{align}\label{e35--1}
	\|\phi_j\|_{H_0^1(\Omega)}^2\le C\|\phi_j\|_{H^{s+1/2-\varepsilon}(\Omega)}^{2(1-\theta)}\|\phi_j\|_{H_0^s(\Omega)}^{2\theta}\le C \|\phi_j\|_{D((-\Delta)_D^s)}^{2(1-\theta)}\|\phi_j\|_{H_0^s(\Omega)}^{2\theta}.
\end{align}
Using the fact that 
\begin{align*}
	\|\phi_j\|_{H_0^1(\Omega)}^2\simeq \int_{\Omega}|\nabla \phi_j|^2\,dx \quad\mbox{ and }\quad  \|\phi_j\|_{H_0^s(\Omega)}^2\simeq \int_{\Omega}|\nabla^s \phi_j|^2\,dx=\lambda_j,
\end{align*}
we can then deduce from \eqref{e35--1} that
\begin{align*}
	\int_{\Omega}|\nabla \phi_j|^2\,dx &\le C\left(\int_{\Omega}|\nabla^s \phi_j|^2\,dx\right)^{\theta}\|\phi_j\|_{D((-\Delta)_D^s)}^{2(1-\theta)} 
	\\
	&= C\left(\int_{\Omega}|\nabla^s \phi_j|^2\,dx\right)\left(\int_{\Omega}|\nabla^s \phi_j|^2\,dx\right)^{\theta-1}\|\phi_j\|_{D((-\Delta)_D^s)}^{2(1-\theta)} 
	\\
	&= C\lambda_j^{2(1-\theta)+\theta-1}\left(\int_{\Omega}|\nabla^s \phi_j|^2\,dx\right) = C\lambda_j^{1-\theta} \int_{\Omega}|\nabla^s \phi_j|^2\,dx 
	\\
	&= C\lambda_j^{\frac{2(1-s)}{1-2\varepsilon}} \int_{\Omega}|\nabla^s \phi_j|^2\,dx,
\end{align*}
and the proof is finished.
\end{proof}

\begin{proposition}\label{propInterp}
Let $\Omega\subset\RR^N$ ($N\ge 1$) be a bounded smooth open set and $s\in[1/4,1/2)$. Then, there is a constant $C=C(N,s,\Omega)>0$ such that for every $j\in\NN$ the eigenfunction $\phi_j$ satisfies
\begin{align}\label{eq:interp}
	\|\phi_j\|_{H^{1-2s}(\Omega)}^2\le C\lambda_j^{\alpha(s)}\int_{\Omega}|\nabla^s\phi_j|^2\,dx,
\end{align}
with
\begin{align}\label{eq:alpha}
	\alpha(s)\coloneqq 
	\begin{cases} 
		\displaystyle\frac 1s-3, & \displaystyle \text{if } s\in\left[\frac 14,\frac 13\right)
		\\[15pt]
		\displaystyle 1-\frac{1}{3s}, & \displaystyle \text{if } s\in\left[\frac 13,\frac 12\right).
	\end{cases}
\end{align}
\end{proposition}

\begin{proof}
The proof is similar to the one of Proposition \ref{lem31} and is based on an interpolation argument. Let us first notice that, since $s\in[1/4,1/2)$, we have from Proposition \ref{grubb}(a) that
\begin{align*}
	\phi_j\in D((-\Delta)^s_D)=H_0^{2s}(\Omega)
\end{align*}
We now distinguish two cases.

\medskip\noindent\textbf{Case 1:} $s\in[1/4,1/3)$: define the parameter 
\begin{align*}
	\theta\coloneqq \frac{4s-1}{s}=4-\frac 1s\in[0,1),
\end{align*}
for which we have 
\begin{align*}
	1-2s = \theta s + (1-\theta)2s.
\end{align*}

Using the interpolation between the Hilbert spaces $H^{2s}_0(\Omega)$, $H^{1-2s}(\Omega)$ and $H_0^s(\Omega)$ (see e.g. \cite[Theorem 11.6]{lions1968problemes}), we get that there is a constant $C=C(N,s,\Omega)>0$ such that for every $j\in\NN$ 
\begin{align*}
	\|\phi_j\|_{H^{1-2s}(\Omega)}^2\le C\|\phi_j\|_{H^{2s}_0(\Omega)}^{2(1-\theta)}\|\phi_j\|_{H_0^s(\Omega)}^{2\theta}= C \|\phi_j\|_{D((-\Delta)_D^s)}^{2(1-\theta)}\|\phi_j\|_{H_0^s(\Omega)}^{2\theta}.
\end{align*}
Arguing as in Proposition \ref{lem31}, we obtain that
\begin{align*}
	\|\phi_j\|_{H^{1-2s}(\Omega)}^2 \leq C\lambda_j^\theta\lambda_j^{2(1-\theta)} = C\lambda_j\lambda_j^{2(1-\theta)+\theta-1} = C\lambda_j^{1-\theta} \int_{\Omega}|\nabla^s \phi_j|^2\,dx= C\lambda_j^{\frac 1s -3} \int_{\Omega}|\nabla^s \phi_j|^2\,dx.
\end{align*}

\medskip\noindent\textbf{Case 2:} $s\in[1/3,1/2)$: define the parameter 
\begin{align*}
	\theta\coloneqq \frac{3s-1}{3s}=1-\frac{1}{3s}\in[0,1),
\end{align*}
for which we have 
\begin{align*}
	1-2s = (1-\theta)s - 2\theta s.
\end{align*}

Using the interpolation between the Hilbert spaces $H^{-2s}(\Omega)$, $H^{1-2s}(\Omega)$ and $H_0^s(\Omega)$ (see e.g. \cite[Theorem 12.3]{lions1968problemes}), together with the fact that $D((-\Delta)^s_D)\hookrightarrow H^{-2s}(\Omega)$, we get that there is a constant $C=C(N,s,\Omega)>0$ such that for every $j\in\NN$ 
\begin{align*}
	\|\phi_j\|_{H^{1-2s}(\Omega)}^2\le C\|\phi_j\|_{H^s_0(\Omega)}^{2(1-\theta)}\|\phi_j\|_{H^{-2s}(\Omega)}^{2\theta}\leq C\|\phi_j\|_{H^s_0(\Omega)}^{2(1-\theta)}\|\phi_j\|_{D((-\Delta)^s_D)}^{2\theta}.
\end{align*}
Arguing as in Proposition \ref{lem31}, we obtain that
\begin{align*}
	\|\phi_j\|_{H^{1-2s}(\Omega)}^2 \leq C\lambda_j^{1-\theta}\lambda_j^{2\theta} = C\lambda_j\lambda_j^\theta = C\lambda_j^\theta \int_{\Omega}|\nabla^s \phi_j|^2\,dx= C\lambda_j^{1-\frac{1}{3s}} \int_{\Omega}|\nabla^s \phi_j|^2\,dx,
\end{align*}
and the proof is finished.
\end{proof}

\begin{proposition}\label{prop:L-infty}
Let $\Omega\subset\RR^N$ ($N\ge 1$) be a bounded smooth open set and $s\in(0,1/4)$. Then, there is a constant $C=C(N,s,\Omega)>0$ such that for every $j\in\NN$ the eigenfunction $\phi_j$  satisfies
\begin{align}\label{eq:L-infty}
	\|\phi_j\|_{L^\infty(\Omega)}\le C\lambda_j^{\frac{N}{4s}}.
\end{align}
\end{proposition}

\begin{proof}
The proof uses a bootstrap argument based on  Proposition \ref{prop:DP_Lebesgue} (see also \cite[Proposition 1.4]{ros2014extremal}). First of all, we notice that for $s\in(0,1/4)$ we have that
\begin{align}
	\frac{N}{2s}>2N\geq 2.
\end{align}
Hence, Proposition \ref{prop:DP_Lebesgue}(b) yields
\begin{align*}
	\|\phi_j\|_{L^{\frac{2N}{N-2s}}(\Omega)}\leq C(N,s)\lambda_j\|\phi_j\|_{L^2(\Omega)} = C(N,s)\lambda_j.
\end{align*}
Now, set $q_0=2$ and define recursively
\begin{align*}
	q_{k+1} = \frac{N q_k}{N-2s q_k}.	
\end{align*}
Let $k^\ast\coloneqq \max\{k\in\NN\,:\,q_k<N/(2s)\}$. By applying iteratively Proposition \ref{prop:DP_Lebesgue}(b), we get
\begin{align*}
	\|\phi_j\|_{L^{q_{k^\ast+1}}(\Omega)} = C(N,s)\lambda_j^{k^\ast}.
\end{align*}
Moreover, since $k^\ast+1>N/(2s)$, we can apply Proposition \ref{prop:DP_Lebesgue}(c) to conclude that
\begin{align}\label{eq:L-infty_prel}
	\|\phi_j\|_{L^\infty(\Omega)}\leq C(N,s,\Omega)\lambda_j\|\phi_j\|_{L^{q_{k^\ast+1}}(\Omega)} \leq C(N,s,\Omega)\lambda_j^{k^\ast+1}.
\end{align}

To conclude the proof, let us analyze the behavior of the sequence $\{q_k\}$ so to determine the value of $k^\ast$. To this end, for all $k\in\NN$ let us denote $\zeta_k = q_k^{-1}$, so that
\begin{align*}
	\zeta_{k+1} = q_{k+1}^{-1} = \frac{1-\gamma q_k}{q_k} = q_k^{-1}-\gamma = \zeta_k - \gamma, \quad \text{with } \gamma := \frac{2s}{N}.	
\end{align*}
Hence, for all $k\in\NN$,
\begin{align*}
	\zeta_k = \zeta_0-k\gamma=\frac 12-k\gamma\quad\Rightarrow\quad q_k=\frac{2}{1-2k\gamma}.
\end{align*}
In particular, we can find the maximal integer $k^\ast$ such that $q_{k^\ast} < N/(2s)$ by solving the inequality
\begin{align*}
	\frac{2}{1 - 2k^\ast\gamma} < \frac{1}{\gamma} \quad \Leftrightarrow \quad k^\ast < \frac{1}{2\gamma}-1 = \frac{N}{4s}-1,	
\end{align*}
that is,
\begin{align*}
	k^\ast = \left\lfloor \frac{N}{4s}-1 \right\rfloor.	
\end{align*}
Using this in \eqref{eq:L-infty_prel}, we finally get \eqref{eq:L-infty}.
\end{proof}

\begin{proposition}\label{prop:L-1_gradient}
Let $\Omega\subset\RR^N$ ($N\ge 1$) be a bounded smooth open set and $s\in(0,1/4)$. Then, there is a constant $C=C(N,s,\Omega)>0$ such that for every $j\in\NN$ the eigenfunction $\phi_j$  satisfies
\begin{align}\label{eq:L1_gradient}
	\|\nabla\phi_j\|_{L^1(\Omega)}\le C\lambda_j^{\frac{N}{4s}+1}.
\end{align}
\end{proposition}

\begin{proof}
Let us start by recalling that, thanks to Proposition \ref{prop:DP_Holder_boundary} (see also \cite[Theorem 1.2]{RS1}), we know that 
\begin{align*}
	v\coloneqq\frac{\phi_j}{\rho^s}\in C^\alpha(\overline\Omega)\;\text{ for  } \alpha<\min\{s,1-s\}, \quad\text{ and } \quad\|v\|_{C^\alpha(\overline\Omega)}\leq C(s,\Omega)\lambda_j\|\phi_j\|_{L^\infty(\Omega)},
\end{align*}
where $\rho(x)=\mbox{dist}(x,\partial\Omega)$ for all $x\in\Omega$. Now, let us write $\phi_j = \rho^sv$, so that
\begin{align*}
	\nabla\phi_j = \nabla(\rho^s v) = \rho^s\nabla v + sv\rho^{s-1}\nabla\rho.
\end{align*}

Using that $\alpha\in (0,1)$ and $v\in {C^\alpha(\overline\Omega)}$ we have that $|\nabla v|\leq C(s,\Omega)\rho^{\alpha-1}\|v\|_{C^\alpha(\overline\Omega)}$. So, we can estimate
\begin{align*}
	|\rho^s\nabla v| = \rho^s|\nabla v| \leq \rho^{s+\alpha-1}\quad\|v\|_{C^\alpha(\overline\Omega)}\leq C(s,\Omega)\rho^{s-1}\lambda_j\|\phi_j\|_{L^\infty(\Omega)}.
\end{align*}
Moreover, since $\rho$ is a Lipschitz continuous function, we also have that
\begin{align*}
	|sv\rho^{s-1}\nabla\rho| = s|v|\rho^{s-1}|\nabla\rho|\leq C(s,\Omega)|v|\rho^{s-1}\leq C(s,\Omega)\rho^{s-1}\lambda_j\|\phi_j\|_{L^\infty(\Omega)}.
\end{align*}
Hence, 
\begin{align*}
	|\nabla\phi_j| \leq C(s,\Omega)\rho^{s-1}\lambda_j\|\phi_j\|_{L^\infty(\Omega)},
\end{align*}
from which we can deduce that 
\begin{align*}
	\|\nabla\phi_j\|_{L^1(\Omega)} = \int_\Omega |\nabla\phi_j|\,dx \leq C(s,\Omega)\lambda_j\|\phi_j\|_{L^\infty(\Omega)}\int_\Omega \rho^{s-1}\,dx \leq C(s,\Omega)\lambda_j\|\phi_j\|_{L^\infty(\Omega)}.
\end{align*}
Using \eqref{eq:L-infty} in this last estimate, we finally get \eqref{eq:L1_gradient}.
\end{proof}

\section{Technical proofs}\label{appendix_proofs}

\begin{proof}[\bf Proof of Proposition \ref{lemmaMult}]
We prove  the identity \eqref{e210} in three steps, as follows.
\begin{itemize}
	\item[1.] Firstly, we obtain a Pohozaev-type identity for the solutions of the fractional wave equation \eqref{adjointWave}.  This is done formally, by employing the Pohozaev identity for the fractional Laplace operator (for the stationary problem) obtained in \cite[Proposition 1.6]{Ro-Se}.
	\item[2.] Secondly, we give a rigorous justification to the computations in the first step, exploiting the fact that solutions $p\in\mathcal H_J$ involve a finite number of Fourier components. The key point here is that $p\in \mathcal H_J$ is a finite spectral packet. Therefore all spatial derivatives and boundary traces involved below inherit the regularity of the eigenfunctions, allowing us to justify rigorously every step of the multiplier argument.
	\item[3.] Finally, we employ a standard argument regarding the equipartition of the energy to obtain our identity \eqref{e210}.
\end{itemize}

\noindent\textbf{Step 1: The Pohozaev identity for solutions of \eqref{adjointWave}}.  The first step is to obtain  the identity
\begin{align}\label{e23}
	\frac{\Gamma(1+s)^2}{2}\int_0^T\int_{\partial\Omega}(x\cdot\nu)\left|\frac{p}{\rho^s}\right|^2\,d\sigma dt =&\, \frac{N}{2}\int_0^T\int_{\Omega}|p_t|^2\,dxdt + \frac{2s-N}{2}\int_0^T\int_{\Omega}|\nabla^sp|^2\,dxdt \notag
	\\
	&+ \int_{\Omega}p_t\big(x\cdot\nabla p\big)\,dx\,\bigg|_{t=0}^{t=T}.
\end{align}

This is done, so far formally, by multiplying the first equation in \eqref{adjointWave} by $x\cdot\nabla p$. Integrating over $\Omega\times(0,T)$ we get that
\begin{align}\label{A1}
	\int_0^T\int_{\Omega}p_{tt}\left(x\cdot\nabla p\right)\,dxdt +\int_0^T\int_{\Omega}(-\Delta)^sp\left(x\cdot\nabla p\right)\,dxdt=0.
\end{align}

Under the regularity assumptions of \cite[Proposition 1.6]{Ro-Se}, that, as we shall see, are fulfilled by solutions of the fractional wave equation involving a finite number of Fourier components, we have that, for all $t\in [0, T]$, 
\begin{align}\label{A2}
	\int_{\Omega}(-\Delta)^sp \left(x\cdot\nabla p\right)\,dx = \frac{2s-N}{2}\int_{\Omega}|\nabla^sp|^2\,dx - \frac{\Gamma(1+s)^2}{2}\int_{\partial\Omega}(x\cdot\nu)\left|\frac{p}{\rho^s}\right|^2\,d\sigma.
\end{align}
Using \eqref{A2} we get from \eqref{A1} that 
\begin{align}\label{A3}
	\int_0^T\int_{\Omega}p_{tt}\left(x\cdot\nabla p\right)\,dxdt = \frac{N-2s}{2}\int_0^T\int_{\Omega}|\nabla^sp|^2\,dxdt + \frac{\Gamma(1+s)^2}{2}\int_0^T\int_{\partial\Omega}(x\cdot\nu)\left|\frac{p}{\rho^s}\right|^2\,d\sigma dt.
\end{align}
Thus, it is sufficient to handle the term
\begin{align*}
	\int_0^T\int_{\Omega}p_{tt}\left(x\cdot\nabla p\right)\,dxdt.
\end{align*}
Formally, this is done as follows:
\begin{align}\label{A4}
	\int_0^T\!\int_{\Omega}p_{tt}\left(x\cdot\nabla p\right)\,dxdt = &\, \int_{\Omega}p_t\left(x\cdot\nabla p\right)\,dx\,\bigg|_{t=0}^{t=T} 
	- \int_0^T\!\int_{\Omega} p_t\left(x\cdot\nabla p_t\right)\,dxdt  \notag
	\\
	= & \int_{\Omega}p_t\left(x\cdot\nabla p\right)\,dx\,\bigg|_{t=0}^{t=T} 
	-\int_0^T\!\int_{\Omega}x\cdot\nabla\left(\frac{|p_{t}|^2}{2}\right)\,dxdt \notag
	\\
	= & \int_{\Omega} p_t\left(x\cdot\nabla p\right)\,dx\,\bigg|_{t=0}^{t=T} + \frac N2\int_0^T\!\int_{\Omega}|p_t|^2\,dxdt.
\end{align}
Once identities \eqref{A3} and \eqref{A4} are obtained and justified, combining them, we get \eqref{e23}.

\smallskip 
\noindent\textbf{Step 2: Regularity analysis.} The developments of Step 1 are formal. To make them rigorous one needs to justify two fundamental aspects. 
\begin{itemize}
	\item[(a)] First of all, we have to show that the finite energy solutions involving a finite number of Fourier components as defined in \eqref{solFourier} fulfill the required regularity properties for the Pohozaev identity.
	\item[(b)] Secondly, we have to justify the integration by parts \eqref{A4} by showing that all the terms in the identity are well defined.
\end{itemize}

Concerning the first aspect, the required regularity is guaranteed by the results of \cite{RS1,Ro-Se,servadei2014weak}. Indeed, according to \cite[Theorem 1.4]{Ro-Se}, if $f\in C_{\rm loc}^{0,1}(\RR)$, then for all $s\in(0,1)$, every weak solution $u\in H_0^s(\Omega)\cap L^\infty(\Omega)$ of the Dirichlet problem
\begin{align}\label{nle}
	(-\Delta)^su=f(u)\;\;\mbox{ in }\;\Omega,\quad u=0\;\mbox{ in }\;\RR^N\setminus\Omega,
\end{align}
is smooth enough and  satisfies the regularity assumptions in \cite[Proposition 1.6]{Ro-Se} to apply the Pohozaev identity \eqref{A2}. 

This is the case for each eigenfunction $(\phi_j)_{j\in\NN}$. Indeed, the nonlinearity $f$ in that case is actually linear, $f(u)=\lambda u$, the eigenfunctions belong to $H_0^s(\Omega)$ by standard well-posedness results (see Proposition \ref{WPprop}), and  the results of \cite{RS1,servadei2014weak} show that they are also bounded a.e. in $\Omega$ and that they belong to the space of H\"older continuous functions $C^{0,s}(\overline{\Omega})$. 

Since each eigenfunction $\phi_j$ satisfies the regularity assumptions in \cite[Proposition 1.6]{Ro-Se} needed to apply \eqref{A2}, the same follows for all finite linear combinations of them, and therefore for the solutions $p(\cdot,t)\in\mathcal H_J$ of the fractional wave equation under consideration.

Now that we know that solutions $p(\cdot,t)\in\mathcal H_J$ of the fractional wave equation are regular enough to apply the Pohozaev identity, let us conclude our discussion by showing that all the terms in \eqref{A4} are well-defined. Firstly, since $p\in C([0,T];H_0^s(\Omega))\cap C^1([0,T];L^2(\Omega))$, we immediately have that
\begin{align*}
	\int_0^T\int_{\Omega} |p_t|^2\,dxdt < +\infty.
\end{align*}
As for the other two integrals in \eqref{A4}, namely
\begin{align}\label{int2}
	\int_0^T\int_\Omega p_{tt}\left(x\cdot\nabla p\right)\,dxdt\quad\mbox{ and }  \int_\Omega p_t\left(x\cdot\nabla p\right)\,dx\,\bigg|_{t=0}^{t=T},
\end{align}
their convergence is a slightly more delicate issue. 

Let us recall that we are considering solutions $p(\cdot,t)\in \mathcal H_J$ which, for a.e. $t\in [0,T]$, are as regular as the eigenfunctions $(\phi_j)_{j\in\NN}$. In particular, since $\phi_j\in L^\infty(\Omega)$ for all $j\in\NN$, we also have that for a.e. $t\in [0,T]$
\begin{align}\label{pLinfty}
	p_t(\cdot,,t) = \sum_{j=1}^J p_j'(t)\phi_j\in L^\infty(\Omega) \quad\text{ and }\quad p_{tt}(\cdot,t) = \sum_{j=1}^J p_j''(t)\phi_j\in L^\infty(\Omega).
\end{align}

Furthermore, since for all $j\in\NN$ the eigenfunction $\phi_j$ is a solution of the elliptic problem \eqref{nle} with $f(\phi_j)=\lambda_j\phi_j$, it follows from \cite[Theorem 1.1]{RSV-CPDE} that 
\begin{align*}
	\frac{\phi_j}{\rho^s}\in C^\gamma(\overline{\Omega}) \text{ for all } \gamma\in (0,s) \quad\text{ and }\quad |\nabla \phi_j|\le C\rho^{s-1}\;\mbox{ in }\Omega, 
\end{align*} 
where $C=C(N,s,\Omega,j)>0$ is a positive constant which depends on $N$, $s$, $\Omega$ and $\phi_j$ (hence, on the frequency $j$). Consequently, 
\begin{align*}
	\frac{p}{\rho^s}\in L^\infty((0, T); C^\gamma(\overline{\Omega}))\text{ for all }\gamma\in(0,s)
\end{align*}
and
\begin{align}\label{pW11}
	\int_\Omega|\nabla \phi_j|\,dx\le C(N,s,\Omega,j)\int_\Omega\rho^{s-1}\,dx<+\infty.
\end{align}

Using \eqref{pLinfty} and \eqref{pW11} we then have that there is a constant $C(\Omega)>0$ such that for a.e. $t\in [0,T]$,
\begin{align*}
	&\left|\int_\Omega p_t(x\cdot\nabla p)\,dx\,\right| \leq C(\Omega) \int_\Omega |p_t|\,|\nabla p|\,dx \leq C(\Omega)  \|p_t(\cdot,t)\|_{L^\infty(\Omega)}\|\nabla p(\cdot,t)\|_{L^1(\Omega)}\leq C(N,s,\Omega,J) 
	\\[8pt]
	&\left|\int_\Omega p_{tt}(x\cdot\nabla p)\,dx\,\right| \leq C(\Omega) \int_\Omega |p_{tt}|\,|\nabla p|\,dx \leq C(\Omega) \|p_{tt}(\cdot,t)\|_{L^\infty(\Omega)}\|\nabla p(\cdot,t)\|_{L^1(\Omega)}\leq C(N,s,\Omega,J). 
\end{align*}
Finally, we can conclude that
\begin{align*}
	\left|\int_0^T\int_\Omega p_{tt}\left(x\cdot\nabla p\right)\,dxdt\,\right| \leq C(N,s,\Omega,J, T) 
\end{align*}
and
\begin{align*}
	\left|\int_\Omega p_t\left(x\cdot\nabla p\right)\,dx\,\bigg|_{t=0}^{t=T}\,\right| \leq C(N,s,\Omega,J, T).
\end{align*}

Hence, the two integrals in \eqref{int2} are all well-defined. Nevertheless, we have to stress that the above estimates are true for solutions $p\in \mathcal H_J$ involving a finite number of Fourier components, and are frequency-dependent. Whether the same can be obtained for general finite energy solutions of \eqref{adjointWave} is an open question, but it is unlikely to be the case due to intervention of higher order terms involving, for instance, $\nabla p$.

\smallskip 
\noindent\textbf{Step 3: Equipartition of the energy}. We now derive the equipartition of the energy identity. Multiplying the first equation in \eqref{adjointWave} by $p$ and integrating by parts over $\Omega\times (0,T)$ (by using the integration by parts formula \eqref{IPF}) we obtain 
\begin{align}\label{e28}
	-\int_0^T\int_{\Omega}|p_t|^2\,dxdt + \int_0^T\int_{\Omega}|\nabla^sp|^2\,dxdt + \int_{\Omega} p_tp\,dx\,\bigg|_{t=0}^{t=T}=0.
\end{align}
Moreover, the conservation of energy yields
\begin{align*}
	\frac N2 & \int_0^T \int_{\Omega}|p_t|^2\,dxdt + \frac{2s-N}{2}\int_0^T\int_{\Omega} |\nabla^sp|^2\,dxdt 
	\\
	&= \frac s2\int_0^T\int_{\Omega}|p_t|^2\,dxdt + \frac s2\int_0^T\int_{\Omega} |\nabla^sp|^2\,dxdt + \frac{N-s}{2}\int_0^T\int_{\Omega}\Big(|p_t|^2 - |\nabla^sp|^2\Big)\,dxdt
	\\
	&= sTE_s(0) + \frac{N-s}{2}\int_0^T\int_{\Omega}\Big(|p_t|^2 - |\nabla^sp|^2\Big)\,dxdt.
\end{align*}
Using this in the identity \eqref{e23}, we then have
\begin{align}\label{e29}
	\frac{\Gamma(1+s)^2}{2}\int_0^T\int_{\partial\Omega}(x\cdot\nu)\left|\frac{p}{\rho^s}\right|^2\,d\sigma dt =&\, sTE_s(0)+ \frac{N-s}{2}\int_0^T\int_{\Omega}\Big(|p_t|^2-|\nabla^sp|^2\Big)\,dxdt \notag
	\\
	& + \int_{\Omega} p_t\left(x\cdot\nabla p\right)\,dx\,\bigg|_{t=0}^{t=T}.
\end{align}
Combining \eqref{e28} and \eqref{e29} we get the identity \eqref{e210} and the proof is finished.
\end{proof}

\begin{proof}[\bf Proof of Proposition \ref{prop31}]
Throughout the proof, the value of the constant $C$ may change from line to line. More precisely, in Case 1 below, corresponding to $s\in(1/2,1)$, we shall write $C=C(N,s,\Omega,\varepsilon)>0$, while in the remaining cases one has $C=C(N,s,\Omega)>0$. Moreover, following the discussion in Remark \ref{rem22}, we address separately the cases $s\in(1/2,1)$, $s=1/2$, $s\in[1/4,1/2)$, and $s\in(0,1/4)$. 

\smallskip 
\noindent\textbf{Case 1:} $s\in(1/2,1)$. From \eqref{e210} we have that
\begin{align}\label{mm}
	sTE_s(0) =\frac{\Gamma(1+s)^2}{2}\int_0^T\int_{\partial\Omega}(x\cdot\nu)\left|\frac{p}{\rho^s}\right|^2\,d\sigma dt 
	- \int_{\Omega}p_t\left(x\cdot\nabla p + \frac{N-s}{2}p\right)\,dx\,\bigg|_{t=0}^{t=T} .
\end{align}

In this case, the dependence on $\varepsilon$ comes from the interpolation estimate \eqref{e35}, and therefore all the constants below are of the form $C=C(N,s,\Omega,\varepsilon)>0$. 

Since we are considering solutions $p\in\mathcal H_J$, and we know from Proposition \ref{grubb}(c) that  $\phi_j\in H_0^1(\Omega)$ for $s\in(1/2,1)$, all terms in \eqref{mm}, and in particular 
\begin{align*}
	\xi(t):= \int_{\Omega}p_t\left(x\cdot\nabla p + \frac{N-s}{2}p\right)\,dx,
\end{align*} 
are well defined. Moreover, by means of Young's inequality, for any $a>0$ we get
\begin{align*}
	|\xi(t)|\le \frac{R}{2a} \int_\Omega |p_t|^2\,dx + \frac{a}{2R} \int_\Omega \left|x\cdot\nabla p + \frac{N-s}{2}p\,\right|^2\,dx,
\end{align*}
where we have denoted 
\begin{align}\label{R_def}
	R\coloneqq \max_{x\in\overline\Omega} |x|. 
\end{align}
Furthermore
\begin{align}\label{e48}
	\int_\Omega \bigg|x\cdot\nabla p + \frac{N-s}{2}p\,\bigg|^2\,dx & = \int_\Omega \left|x\cdot\nabla p\,\right|^2\,dx + \left(\frac{N-s}{2}\right)^2\int_\Omega |p|^2\,dx + (N-s)\int_\Omega p(x\cdot\nabla p)\,dx \notag 
	\\
	& = \int_\Omega \left|x\cdot\nabla p\,\right|^2\,dx + \left(\frac{N-s}{2}\right)^2\int_\Omega |p|^2\,dx + \frac{N-s}{2}\int_\Omega x\cdot\nabla(p^2)\,dx \notag
	\\
	& = \int_\Omega \left|x\cdot\nabla p\,\right|^2\,dx + \left(\frac{N-s}{2}\right)^2\int_\Omega |p|^2\,dx - \frac{N(N-s)}{2}\int_\Omega |p|^2\,dx\notag 
	\\
	& = \int_\Omega \left|x\cdot\nabla p\,\right|^2\,dx - \frac{(N^2-s^2)}{4}\int_\Omega |p|^2\,dx.
\end{align}
Since $N\geq 1$ and $0<s<1$, the last term in the right hand side of \eqref{e48} is negative. Hence, 
\begin{align*}
	\int_\Omega \left|x\cdot\nabla p + \frac{N-2s}{2}p\,\right|^2\,dx \leq \int_\Omega \left|x\cdot\nabla p\,\right|^2\,dx \leq R^2\int_\Omega \left|\nabla p\,\right|^2\,dx.
\end{align*}
This yields
\begin{align}\label{uu-b}
	|\xi(t)|\le \frac{R}{2a} \int_\Omega |p_t|^2\,dx + \frac{aR}{2} \int_\Omega |\nabla p|^2\,dx.
\end{align}
It follows from \eqref{e35} that there is a constant $C=C(N,s,\Omega,\varepsilon)>0$ such that for all $ p\in \mathcal H_J$,
\begin{align}\label{e312}
	\int_\Omega |\nabla p|^2\,dx \leq C\lambda_J^{\frac{2(1-s)}{1-2\varepsilon}}\int_\Omega |\nabla^s p|^2\,dx.
\end{align}
Combining \eqref{uu-b}-\eqref{e312} we then obtain that
\begin{align*}
	\left|\xi(t)\right|\leq \left(\frac{R}{2a} + C aR \lambda_J^{\frac{2(1-s)}{1-2\varepsilon}}\right) E_s(0).
\end{align*}
Choosing now $a = \lambda_J^{-\frac{1-s}{1-2\varepsilon}}$, and since $\Omega$ is bounded, we obtain from the previous estimate that
\begin{align}\label{est_xi2}
	\left|\xi(t)\right|\leq C \lambda_J^{\frac{1-s}{1-2\varepsilon}} E_s(0).
\end{align}
We employ \eqref{est_xi2} to bound the last  term on the right-hand side of \eqref{mm}. In particular, we have 
\begin{align*}
	\bigg| \int_{\Omega} p_t \bigg(x \cdot\nabla p + \frac{N-s}{2}p\bigg)\,dx\,\bigg|_{t=0}^{t=T}\,\bigg|= \bigg|\xi(t)\,\bigg|_{t=0}^{t=T}\bigg| \leq  2C\lambda_J^{\frac{1-s}{1-2\varepsilon}} E_s(0).
\end{align*}
Consequently, from \eqref{mm} we obtain that
\begin{align*}
	\Big(sT - C\lambda_J^{\frac{1-s}{1-2\varepsilon}}\Big)E_s(0) & \le \frac{\Gamma(1+s)^2}{2}\int_0^T\int_{\partial\Omega}(x\cdot\nu)\left|\frac{p}{\rho^s}\right|^2\,d\sigma dt
	\\
	&= \frac{\Gamma(1+s)^2}{2}\left(\int_0^T\int_{\partial\Omega^+}(x\cdot\nu)\left|\frac{p}{\rho^s}\right|^2\,d\sigma dt + \int_0^T\int_{\partial\Omega^-}(x\cdot\nu)\left|\frac{p}{\rho^s}\right|^2\,d\sigma dt\right)
	\\
	&\leq \frac{\Gamma(1+s)^2}{2}\int_0^T\int_{\partial\Omega^+}(x\cdot\nu)\left|\frac{p}{\rho^s}\right|^2\,d\sigma dt.
\end{align*}
Thus, if $T>T_0(J)$ with $T_0(J)$ given by \eqref{e39}, we finally obtain
\begin{align*}
	E_s(0) \leq \frac{\Gamma(1+s)^2}{2s(T-T_0(J))}\int_0^T\int_{\partial\Omega^+}(x\cdot\nu)\left|\frac{p}{\rho^s}\right|^2\,d\sigma dt.
\end{align*}

\noindent\textbf{Case 2:} $s\in(0,1/2]$. As anticipated in Remark \ref{rem22}, when $s\in(0,1/2]$ we shall view the term
\begin{align*}
	\int_{\Omega} p_t\left(x\cdot\nabla p + \frac{N-s}{2}p\right)\,dx
\end{align*} 
as the duality map
\begin{align*}
	\left\langle p_t, x\cdot\nabla p + \frac{N-s}{2}p\right\rangle_{V, V^\star},
\end{align*} 
where $V:=D((-\Delta)_D^s)$ and $V^\star$ is the dual of $V$ with respect to the pivot space $L^2(\Omega)$. When doing that, \eqref{e210} becomes
\begin{align}\label{mm-bis}
	sTE_s(0) =&\, \frac{\Gamma(1+s)^2}{2}\int_0^T\int_{\partial\Omega}(x\cdot\nu)\left|\frac{p}{\rho^s}\right|^2\,d\sigma dt 
	- \Bigg\langle p_t,\left(x\cdot\nabla p+ \frac{N-s}{2}p\right)\Bigg\rangle_{V, V^\star}\,\Bigg|_{t=0}^{t=T} .
\end{align}
Define
\begin{align*}
	\xi(t):=\Bigg\langle p_t,x\cdot\nabla p+ \frac{N-s}{2}p\Bigg\rangle_{V, V^\star} =\underbrace{\langle p_t,x\cdot\nabla p\rangle_{V, V^\star}}_{\xi_1(t)} + \underbrace{\Bigg( p_t,\frac{N-s}{2}p\Bigg)_{L^2(\Omega)}}_{\xi_2(t)}.
\end{align*} 
Then, 
\begin{align}\label{xi12}
	|\xi(t)|\leq |\xi_1(t)|+|\xi_2(t)|.
\end{align}
Let us now estimate the two terms $|\xi_1(t)|$ and $|\xi_2(t)|$. For the latter one, Young's inequality easily yields
\begin{align}\label{xi2_est}
	|\xi_2(t)|\le \frac 12 \|p_t\|_{L^2(\Omega)}^2 + C(N,s,\Omega)\|p\|_{L^2(\Omega)}^2\leq C(N,s,\Omega)E_s(0).
\end{align}
Concerning $|\xi_1(t)|$, instead, we have that for any $a>0$
\begin{align}\label{xi1_est}
	|\xi_1(t)|\le \frac{R}{2a}\|p_t\|_V^2 + \frac {a}{2R}\|x\cdot\nabla p\|_{V^\star}^2 \le \frac{R}{2a}\|p_t\|_V^2 + \frac {aR}{2}\|\nabla p\|_{V^\star}^2,
\end{align}
with $R$ defined in \eqref{R_def}. We now have to consider the cases $s=1/2$, $s\in[1/4,1/2]$, and $s\in(0,1/4)$ separately. 

\noindent\textbf{Sub}-\textbf{case 2.1:} $s= 1/2$. Using \eqref{domain_char}, \eqref{solFourier}, the orthonormality of the eigenfunctions in $L^2(\Omega)$ and the Bessel inequality, we get that 
\begin{align}\label{pt_norm}
	\|p_t\|_V^2 = \sum_{j=1}^J |\lambda_j (p_t,\phi_j)_{L^2(\Omega)}|^2 \leq \lambda_J^2\|p_t\|_{L^2(\Omega)}^2.
\end{align}
Now, by Proposition \ref{grubb}(b), for every $\varepsilon>0$ we have that
\begin{align*}
	H^{-2s+\varepsilon}(\Omega)=H^{-1+\varepsilon}(\Omega)\hookrightarrow V^\star.	
\end{align*}
Hence, 
\begin{align}\label{nabla_p_norm_prel_12}
	\|\nabla p\|_{V^\star}^2 \leq C\|\nabla p\|_{H^{-1+\varepsilon}(\Omega)}^2 \leq C(N,\Omega)\|p\|_{H^{\varepsilon}(\Omega)}^2. 
\end{align}
Since the above estimate holds for any $\varepsilon>0$, we can in particular fix $\varepsilon = 1/2$ so that from \eqref{nabla_p_norm_prel_12} we get
\begin{align}\label{grad_norm12}
	\|\nabla p\|_{V^\star}^2 \leq C(N,\Omega)\|p\|_{H^{\frac 12}(\Omega)}^2 \leq C(N,\Omega)\|\nabla^{\frac 12} p\|_{L^2(\Omega)}^2.
\end{align}
In view of this, using \eqref{xi1_est}, \eqref{pt_norm} and \eqref{grad_norm12}, we obtain 
\begin{align*}
	|\xi_1(t)| \leq C(N,\Omega)\left(\frac{\lambda_J^2}{a} + a \right)E_{\frac 12}(0).
\end{align*}
Taking $a=\lambda_J$ in the above estimate, we then get
\begin{align*}
	|\xi_1(t)| \leq C(N,\Omega)\lambda_J E_{\frac 12}(0)
\end{align*}
which, combined with \eqref{xi12} and \eqref{xi2_est}, yields
\begin{align*}
	|\xi(t)| \leq C(N,\Omega)\lambda_J E_{\frac 12}(0).
\end{align*}
Arguing as in the case $s\in(1/2,1)$, this last inequality implies that 
\begin{align*}
	\bigg[\frac T2-C(N,\Omega)\lambda_J\bigg]E_{\frac 12}(0) \le \frac{\Gamma(3/2)^2}{2}\int_0^T\int_{\partial\Omega^+}(x\cdot\nu)\left|\frac{p}{\rho^{\frac 12}}\right|^2\,d\sigma dt .
\end{align*}
Thus, if $T>T_0(J)$ with $T_0(J)$ defined in \eqref{e39}, we get the estimate \eqref{e37} also for $s=1/2$.  

\noindent\textbf{Sub}-\textbf{case 2.2:} $s\in[1/4,1/2)$. For this range of $s$, according to Proposition \ref{grubb} we have the continuous embedding $V \hookrightarrow H^{2s}_0(\Omega)$ and, consequently, $H^{-2s}(\Omega)\hookrightarrow V^\star$. Hence, we get that
\begin{align}\label{nabla_p_norm_prel}
	\|\nabla p\|_{V^\star}^2 \leq C(N,s,\Omega)\|\nabla p\|_{H^{-2s}(\Omega)}^2 \leq C(N,s,\Omega)\|p\|_{H^{1-2s}(\Omega)}^2.
\end{align}
In addition, using \eqref{solFourier} and \eqref{eq:interp},  we get that there is a constant $C(N,s,\Omega)>0$ such that
\begin{align}\label{bb1}
	\|p\|_{H^{1-2s}(\Omega)}^2 \leq C(N,s,\Omega)\lambda_J^{\alpha(s)} \|\nabla^s p\|_{L^2(\Omega)}^2,
\end{align}
with $\alpha(s)$ given by \eqref{eq:alpha}. Plugging \eqref{bb1} into \eqref{nabla_p_norm_prel} yields
\begin{align}\label{grad_norm}
	\|\nabla p\|_{V^\star}^2 \leq C(N,s,\Omega)\lambda_J^{\alpha(s)} \|\nabla^s p\|_{L^2(\Omega)}^2.
\end{align}
Hence, using \eqref{xi1_est}, \eqref{pt_norm} and \eqref{grad_norm}, we obtain 
\begin{align*}
	|\xi_1(t)| \leq C(N,s,\Omega)\left(\frac{\lambda_J^2}{a} + a\lambda_J^{\alpha(s)} \right)E_s(0).
\end{align*}
Taking $a=\lambda_J^{1-\alpha(s)/2}$ in the above estimate, we then get
\begin{align*}
	|\xi_1(t)| \leq C(N,\Omega)\lambda_J^{1+\frac{\alpha(s)}{2}} E_s(0)
\end{align*}
which, combined with \eqref{xi12} and \eqref{xi2_est}, yields
\begin{align*}
	|\xi(t)| \leq C(N,\Omega)\lambda_J^{1+\frac{\alpha(s)}{2}} E_s(0).
\end{align*}
Arguing as before, this latter inequality implies that 
\begin{align*}
	\bigg[sT-C(N,s,\Omega)\lambda_J^{1+\frac{\alpha(s)}{2}}\bigg]E_s(0) \le \frac{\Gamma(1+s)^2}{2}\int_0^T\int_{\partial\Omega^+}(x\cdot\nu)\left|\frac{p}{\rho^s}\right|^2\,d\sigma dt.
\end{align*}
Thus, if $T>T_0(J)$ with $T_0(J)$ defined in \eqref{e39}, we get once again the desired estimate \eqref{e37}. 

\noindent\textbf{Sub}-\textbf{case 2.3:} $s\in(0,1/4)$. Since $p_t(t,\cdot)\in L^\infty(\Omega)$ and $x\cdot\nabla p(t,\cdot)\in L^1(\Omega)$ (see Step 2 of the proof of Proposition \ref{lemmaMult}), we have that
\begin{align*}
	\xi_1(t) = \langle p_t,x\cdot\nabla p\rangle_{V, V^\star}=\int_{\Omega}p_t(x\cdot\nabla p)\,dx.
\end{align*}
Hence, H\"older and Young's inequalities yield that, for all $a>0$ and $R$ defined in \eqref{R_def},
\begin{align}\label{xi1_est_low}
	|\xi_1(t)|\le \frac{R}{2a}\|p_t\|_{L^\infty(\Omega)}^2 + \frac {aR}{2}\|\nabla p\|_{L^1(\Omega)}^2.
\end{align}
Moreover, using \eqref{solFourier} and the Cauchy-Schwarz and Bessel inequalities, we can estimate
\begin{displaymath}
	\begin{array}{l}
		\displaystyle\|p_t\|_{L^\infty(\Omega)}^2 \leq \left(\sum_{j=1}^J |p_j'(t)|^2\right) \left(\sum_{j=1}^J \|\phi_j\|_{L^\infty(\Omega)}^2\right) \leq C(N,s,\Omega)\|p_t\|_{L^2(\Omega)}^2\left(\sum_{j=1}^J \|\phi_j\|_{L^\infty(\Omega)}^2\right),
		\\[25pt]
		\displaystyle\|\nabla p\|_{L^1(\Omega)}^2 \leq \left(\sum_{j=1}^J |p_j(t)|^2\right) \left(\sum_{j=1}^J \|\nabla\phi_j\|_{L^1(\Omega)}^2\right) \leq C(N,s,\Omega)\|p\|_{L^2(\Omega)}^2\left(\sum_{j=1}^J \|\nabla\phi_j\|_{L^1(\Omega)}^2\right).
	\end{array}
\end{displaymath}
Now, using \eqref{eq:weyl} and \eqref{eq:L-infty}, we can estimate
\begin{align*}
	\sum_{j=1}^J \|\phi_j\|_{L^\infty(\Omega)}^2 \leq C(N,s,\Omega) \sum_{j=1}^J \lambda_j^{\frac{N}{2s}} \leq C(N,s,\Omega) \sum_{j=1}^J j \leq C(N,s,\Omega) J^2 \leq C(N,s,\Omega) \lambda_J^{\frac Ns}.
\end{align*}
In the same way, from \eqref{eq:weyl} and \eqref{eq:L1_gradient}, we have that
\begin{align*}
	\sum_{j=1}^J \|\nabla\phi_j\|_{L^1(\Omega)}^2 \leq C(N,s,\Omega) \sum_{j=1}^J \lambda_j^{\frac{N}{2s}+2} \leq C(N,s,\Omega) \sum_{j=1}^J j^{\frac{4s}{N}+1} \leq C(N,s,\Omega) J^{\frac{4s}{N}+2} \leq C(N,s,\Omega) \lambda_J^{\frac Ns+2}.
\end{align*}
Using this in \eqref{xi1_est_low}, we obtain 
\begin{align*}
	|\xi_1(t)|\le \frac{R}{2a}\lambda_J^{\frac Ns}\|p_t\|_{L^2(\Omega)}^2 + \frac {aR}{2}\lambda_J^{\frac Ns +2}\|p\|_{L^2(\Omega)}^2 \leq C(N,s,\Omega)\left(\frac{\lambda_J^{\frac Ns}}{a} + a\lambda_J^{\frac Ns +2}\right)E_s(0).
\end{align*}
Taking $a=\lambda_J^{-1}$ in the above estimate, we then get
\begin{align*}
	|\xi_1(t)| \leq C(N,\Omega)\lambda_J^{\frac Ns +1} E_s(0)
\end{align*}
which, combined with \eqref{xi12} and \eqref{xi2_est}, yields
\begin{align*}
	|\xi(t)| \leq C(N,\Omega)\lambda_J^{\frac Ns + 1} E_s(0).
\end{align*}
Arguing as before, this latter inequality implies that 
\begin{align*}
	\bigg[sT-C(N,s,\Omega)\lambda_J^{\frac Ns +1}\bigg]E_s(0) \le \frac{\Gamma(1+s)^2}{2}\int_0^T\int_{\partial\Omega^+}(x\cdot\nu)\left|\frac{p}{\rho^s}\right|^2\,d\sigma dt.
\end{align*}
Thus, if $T>T_0(J)$ with $T_0(J)$ defined in \eqref{e39}, we get the estimate \eqref{e37} also for $s\in(0,1/4)$. 
\end{proof}

\begin{proof}[\bf Proof of Proposition \ref{prop62}]
First of all, let us notice that, up to a change of variables $t\mapsto T-t$, the partial observability inequality \eqref{Ex1} is equivalent to the following one:
\begin{align}\label{obs_q}
	\|q(\cdot,T)\|_{L^2(\Omega)}^2\le \frac{\Gamma(1+s)^2}{3sT}\exp\left(\frac{C(s,\Omega)}{T}T_0(J)^2\right)
	\int_0^T\int_{\partial\Omega^+}(x\cdot \nu)\left|\frac{q}{\rho^s}\right|^2\,d\sigma dt,
\end{align}
where $q$ solves the forward equation
\begin{equation}\label{heat_q}
	\begin{cases}
		q_t +(-\Delta)^sq=0 &\mbox{in }\;\Omega\times (0,T),
		\\
		q=0 &\mbox{in }\; (\mathbb R^N\setminus\Omega)\times (0,T),
		\\
		q(\cdot,0)=q_0 &\mbox{in }\;\Omega,
	\end{cases}
\end{equation}
with $q_0\in \mathcal H_J$.

We prove the estimate \eqref{obs_q} by applying the transmutation techniques of \cite{ervedoza1observability,ervedoza2011sharp}. To this end, let $T,L>0$ be given real numbers and let $k:(-L,L)\times(0,T)\to\mathbb{R}$ be the solution of \eqref{eq:k}. Define the function
\begin{align*}
	\psi(x,\zeta) :=\int_0^T k(\zeta,t)q(x,t)\,dt,
\end{align*}
with $q$ the solution of \eqref{heat_q}. Then a simple calculation gives
\begin{align*}
	\psi_{\zeta\zeta} = \int_0^T k_{\zeta\zeta}q\,dt = -\int_0^T k_tq\,dt = \int_0^T kq_t\,dt = -\int_0^Tk(-\Delta)^sq\,dt = -(-\Delta)^s\psi.
\end{align*}
Hence, $\psi$ is a weak solution of the following fractional wave equation:
\begin{align}\label{eq:v}
	\begin{cases}
		\psi_{\zeta\zeta} + (-\Delta)^s\psi = 0 &\mbox{in }\;\Omega\times(-L,L),
		\\
		\psi = 0 &\mbox{in }\;(\mathbb{R}^N\setminus\Omega)\times(-L,L),
		\\
		\displaystyle \psi(x,0) = 0, \quad \psi_\zeta(x,0) = \int_0^T k_\zeta(0,t)q(x,t)\,dt &\mbox{in }\; \Omega.
	\end{cases}
\end{align}

Notice that $\zeta$ plays the role of the space variable in \eqref{eq:k}, while it represents the time variable in \eqref{eq:v}. Now, using \eqref{e37}, if $L> T_0/2$, we have that 
\begin{align}\label{obs_prel}
	\Bigg\|\int_0^T k_\zeta(0,t)q(x,t)\,dt\,\Bigg\|_{L^2(\Omega)}^2 \leq \frac{\Gamma(1+s)^2}{2s(2L-T_0)}\int_{-L}^L\int_{\partial\Omega^+}\int_0^T |k(\zeta,t)|^2\left|\frac{q(x,t)}{\rho(x)^s}\right|^2(x\cdot\nu)\,dtd\sigma d\zeta.
\end{align}
Let us now expand the initial datum $q_0\in \mathcal H_J$ in the basis of the eigenfunctions as follows:
\begin{align*}
	q_0(x) = \sum_{j=1}^J q_j\phi_j(x)\;\;\mbox{ with }\; q_j:=(q_0,\phi_j)_{L^2(\Omega)}.
\end{align*}
Then, the solution of the system \eqref{heat_q} is given by
\begin{align}\label{uu}
	q(x,t) = \sum_{j=1}^J q_je^{-\lambda_j t}\phi_j(x).
\end{align}

Plugging \eqref{uu} in \eqref{obs_prel}, using the orthonormality of the eigenfunctions and \eqref{k_id}, for all $\beta>2L^2$, we obtain the following estimates:
\begin{align*}
	\left\|\int_0^T k_\zeta(0,t)q(x,t)\,dt\,\right\|_{L^2(\Omega)}^2 &= \sum_{j=1}^J |q_j|^2\left[\int_0^T \exp\left(-\beta\left(\frac 1t + \frac{1}{T-t}\right)-\lambda_jt\right)\,dt\right]^2
	\\
	&\geq \sum_{j=1}^J |q_j|^2e^{-2\lambda_jT}\left[\int_{T/4}^{T/2} \exp\left(-\beta\left(\frac 1t + \frac{1}{T-t}\right)\right)\,dt\right]^2
	\\
	&\geq \frac{T^2}{4}\exp\left(-\frac{8\beta}{T}\right)\sum_{j=1}^J |q_j|^2e^{-2\lambda_jT} 
	\\
	&= \frac{T^2}{4}\exp\left(-\frac{8\beta}{T}\right)\|q(\cdot,T)\|_{L^2(\Omega)}^2.
\end{align*}
That is, 
\begin{align}\label{jj}
	\|q(\cdot,T) \|_{L^2(\Omega)}^2 \leq \frac{4}{T^2}\exp\left(\frac{8\beta}{T}\right)\left\|\int_0^T k_\zeta(0,t)q(x,t)\,dt\,\right\|_{L^2(\Omega)}^2. 	
\end{align}
From \eqref{obs_prel} and \eqref{jj} we then get 
\begin{align}\label{obs_prel2}
	\|q(\cdot,T) \|_{L^2(\Omega)}^2 \leq \frac{\Gamma(1+s)^2}{2sT^2(2L-T_0)}\exp\left(\frac{8\beta}{T}\right)\int_{-L}^L\int_{\partial\Omega^+}\int_0^T |k(\zeta,t)|^2\left|\frac{q(x,t)}{\rho(x)^s}\right|^2(x\cdot\nu)\,dtd\sigma d\zeta.
\end{align}
Now, using \eqref{k_est1}, we obtain that
\begin{align*}
	\int_0^T |k(\zeta,t)|^2\left|\frac{q(x,t)}{\rho(x)^s}\right|^2(x\cdot\nu)\,dt \leq |\zeta|^2\left(\int_0^T \!(x\cdot\nu)\left|\frac{q(x,t)}{\rho(x)^s}\right|^2\,dt\right)\left(\int_0^T\!\exp\left[\frac{2}{\min\{t,T-t\}}\left(\frac{\zeta^2}{\delta}-\frac{\beta}{1+\delta}\right)\right]\,dt\right).
\end{align*}
Moreover, since we are taking $\beta>2L^2$, we have that
\begin{align*}
	\frac{L^2}{\beta-L^2} < \frac{L^2}{L^2} = 1,
\end{align*}
which immediately implies that the intervals
\begin{align*}
	\left(\frac{L^2}{\beta-L^2},1\right)\neq\varnothing \quad\text{ and }\quad \left(\frac{L^2}{\beta-L^2},1\right)\subset (0,1).
\end{align*}
In view of this, we can fix
\begin{align*}
	\delta\in \left(\frac{L^2}{\beta-L^2},1\right)
\end{align*} 
so that, for all $\zeta\in (-L,L)$ and $\beta>2L^2$, 
\begin{align*}
	\frac{\zeta^2}{\delta}-\frac{\beta}{1+\delta} &< \zeta^2-\frac{\beta}{1+\frac{L^2}{\beta-L^2}} = \zeta^2-\beta+L^2 < L^2-\beta+L^2 = 2L^2-\beta <0.
\end{align*}
In particular, we get that
\begin{align}\label{obs_prel5}
	\int_0^T |k(\zeta,t)|^2\left|\frac{q(x,t)}{\rho(x)^s}\right|^2(x\cdot\nu)\,dt &\leq |\zeta|^2\left(\int_0^T (x\cdot\nu)\left|\frac{q(x,t)}{\rho(x)^s}\right|^2\,dt\right)\left(\int_0^T \,dt\right) \notag
	\\
	&=|\zeta|^2T\int_0^T (x\cdot\nu)\left|\frac{q(x,t)}{\rho(x)^s}\right|^2\,dt. 
\end{align}
Substituting \eqref{obs_prel5} into \eqref{obs_prel2}, and thanks to Fubini's Theorem, we obtain the following estimates:
\begin{align}\label{obsHeatPrel}
	\|q(\cdot,T) \|_{L^2(\Omega)}^2 &\leq \frac{\Gamma(1+s)^2}{2sT^2(2L-T_0)}\exp\left(\frac{8\beta}{T}\right)\int_{-L}^L\int_{\partial\Omega^+} |\zeta|^2T\int_0^T (x\cdot\nu)\left|\frac{q(x,t)}{\rho(x)^s}\right|^2\,dtd\sigma d\zeta \notag 
	\\
	&= \frac{\Gamma(1+s)^2}{2sT(2L-T_0)}\exp\left(\frac{8\beta}{T}\right) \int_{-L}^L |\zeta|^2\left(\int_0^T\int_{\partial\Omega^+}(x\cdot\nu)\left|\frac{q(x,t)}{\rho(x)^s}\right|^2\,d\sigma dt\right)\,d\zeta \notag 
	\\
	&= \frac{\Gamma(1+s)^2}{3sT(2L-T_0)}L^3\exp\left(\frac{8\beta}{T}\right) \int_0^T\int_{\partial\Omega^+} (x\cdot\nu)\left|\frac{q(x,t)}{\rho(x)^s}\right|^2\,d\sigma dt.
\end{align}
Finally, if we take $L = T_0>T_0/2$ and $\beta=3T_0^2=3L^2>2L^2$, we get from \eqref{obsHeatPrel} that
\begin{align*}
	\|q(\cdot,T) \|_{L^2(\Omega)}^2 &\leq \frac{\Gamma(1+s)^2}{3sTT_0(J)}T_0(J)^3\exp\left(\frac{24}{T}T_0(J)^2\right) \int_0^T\int_{\partial\Omega^+} (x\cdot\nu)\left|\frac{q(x,t)}{\rho(x)^s}\right|^2\,d\sigma dt
	\\
	&= \frac{\Gamma(1+s)^2}{3sT}T_0(J)^2\exp\left(\frac{24}{T}T_0(J)^2\right) \int_0^T\int_{\partial\Omega^+} (x\cdot\nu)\left|\frac{q(x,t)}{\rho(x)^s}\right|^2\,d\sigma dt
	\\	
	&\leq \frac{\Gamma(1+s)^2}{3sT}\exp\left(\frac{C(s,\Omega)}{T}T_0(J)^2\right) \int_0^T\int_{\partial\Omega^+} (x\cdot\nu)\left|\frac{q(x,t)}{\rho(x)^s}\right|^2\,d\sigma dt,
\end{align*}
and the proof is finished.
\end{proof}
}

\bibliographystyle{acm}
\bibliography{biblio}

\end{document}